\documentclass[12pt]{article}
\usepackage{amssymb}
\usepackage{amsmath}
\usepackage{amsthm}
\usepackage{tikz}
\usepackage{amscd}
\usepackage{latexsym}
\usepackage{mathrsfs}
\usepackage{enumitem}
\usepackage{lipsum}
\usepackage[noadjust]{cite}
\usepackage{float}
\usepackage{subfig}
\usepackage{graphicx}
\usepackage[justification =centering, labelsep =period]{caption}

\newtheorem{thm}{Theorem}[section]
\newtheorem{cor}[thm]{Corollary}
\newtheorem{lem}[thm]{Lemma}
\newtheorem{Claim}[thm]{Claim}
\newtheorem{prop}[thm]{Proposition}

\newtheorem{exa}[thm]{Example}

\newtheorem{no}[thm]{Notation}
\newcommand{\pf}{\begin{proof}}
\newcommand{\epf}{\end{proof}}

\numberwithin{equation}{section}

\begin{document}
\title{{\bf The varieties generated by 3-hypergraph semirings} \footnote{This paper is supported by National Natural Science
Foundation of China (11971383)}}
\author{{\bf Yuanfan Zhuo}$^a$, {\bf Xingliang Liang}$^b$, {\bf Yanan Wu}$^a$, {\bf Xianzhong Zhao}$^a$
  \footnote{Corresponding author. E-mail: zhaoxz@nwu.edu.cn} \\
   {\small $^a$ School of Mathematics, Northwest University}\\
   {\small Xi'an, Shaanxi, 710127, P. R. China}\\
   {\small $^b$ School of Mathematics and Data Science},\\
   {\small Shaanxi University of Science and Technology} \\
   {\small Xian, Shaanxi, 710021, P.R. China}\\}
\date{}
\maketitle
\vskip -8pt
\baselineskip 16pt
\newcommand\blfootnote[1]{%
\begingroup
\renewcommand\thefootnote{}\footnote{#1}%
\addtocounter{footnote}{-1}%
\endgroup}
\begin{center}
\begin{minipage}{140mm}
\noindent\textbf{ABSTRACT} In this paper the 3-hypergraph semigroups and 3-hypergraph semirings from 3-hypergraphs $\mathbb{H}$ are introduced and the varieties generated by them are studied. It is shown that all 3-hypergraph semirings $S_{\scriptscriptstyle \mathbb{H}}$ are nonfinitely based and subdirectly irreducible.
Also, it is proved that each variety generated by 3-hypergraph semirings is equal to a variety generated by 3-uniform hypergraph semirings.
It is well known that both variety $\mathbf{V}(S_c(abc))$  (see, J. Algebra 611: 211--245, 2022 and J. Algebra 623: 64--85, 2023) and variety $\mathbf{V}(S_{\scriptscriptstyle \mathbb{H}})$ play key role in the theory of variety of ai-semirings, where 3-uniform hypergraph $\mathbb{H}$ is a 3-cycle.
They are shown that each variety generated by
2-robustly strong 3-colorable 3-uniform hypergraph semirings is equal to variety $\mathbf{V}(S_c(abc))$,
and each variety generated by so-called beam-type hypergraph semirings or fan-type hypergraph semirings is equal to the variety $\mathbf{V}(S_{\scriptscriptstyle \mathbb{H}})$ generated by a 3-uniform 3-cycle hypergraph semiring $S_{\scriptscriptstyle \mathbb{H}}$.
Finally, an infinite ascending chain is provided in the lattice of subvarieties of the variety generated by all 3-uniform hypergraph semirings.
This implies that the variety generated by all 3-uniform hypergraph semirings has infinitely many subvarieties.

\end{minipage}
\end{center}

\section{Introduction}

A \emph{semiring} means an algebra with two associative binary operations $+,\, \cdot,\,$ in which
$+$ is commutative and $\cdot$ distributes over $+$ from the left and right.
A semiring $S$ is called \emph{commutative} if $ab=ba$ for all $a, b\in S$.
An \emph{additively idempotent semiring}(\,ai-semiring for short) is a semiring whose additive reduct is a commutative idempotent semigroup.
Recall that a variety is said to be \emph{finitely based} if it can be defined by finitely many identities; otherwise is said to be \emph{nonfinitely based}.
For any semiring $S$, we shall use $\mathbf{V}(S)$ to denote the variety generated by $S$.
A semiring $S$ is said to be \emph{finitely based} (resp., \emph{nonfinitely based}) if the variety $\mathbf{V}(S)$ is finitely based (resp., nonfinitely based).
Ai-semiring varieties have been frequently studied in the literature (see, \cite{Jackson, pas1, pz2, Ren, rz, rzs1, rzs2, rzw, Shao, Wu3, Wu2, Wu, zxz, zgs, zrp, zxr, zsg}).

By a flat semiring we mean an ai-semiring $S$ such that its multiplicative reduct has a zero element $0$ and $a+b=0$ for all distinct elements $a$ and $b$ of $S$.
Such operation is called the flat addition. In other words, a nontrivial ai-semiring is flat if and only if its additive reduct is a semilattice of height one and the top element is 0. Some progress about flat algebras or flat semirings has been made by Jackson, Zhao, Ren, Wu and Kun (see, \cite{jackson1, Jackson, kun, Ren, Wu2, Wu, zxr}).
Recall that  if $S$ is a semigroup (or semiring) with $0$,
we say that it satisfies the \emph{0-cancellative laws} if for any $ a,\,b,\, c \in S$,
$$ab = ac\neq 0 \Rightarrow b = c, \,\, and \,\, ba = ca\neq 0 \Rightarrow b = c.$$
It is well known that a semigroup with 0 is the multiplicative reduct of a flat semiring if and only if it satisfies
the 0-cancellative laws (see, \cite [Lemma 2.2]{Jackson}).

We know that the two classes of flat semirings play important roles in the theory of semiring. One class of flat semirings is $S_c(W)$ (see, \cite{Jackson}).
Let $X_c^+$ denote the free commutative semigroup on some alphabet $X$, and let $W$ be a nonempty subset of $X_c^+$.
If $S_c(W)$ denotes the set of all nonempty subwords of words in $W$ with an additional symbol $0$ and define an operation on $S_c(W)$ by the rule ($\forall u, v\in S_c(W)$) 
\begin{equation*}
u\cdot v= \begin{cases}
uv, & \text{if} ~ uv \in S_c(W)\setminus \{0\};\\
0, & \text{otherwise}.
\end{cases}
\end{equation*}
Then $S_c(W)$ forms a semigroup with a zero element $0$.
It is easy to see that the semigroup $S_c(W)$ is 0-cancellative.
So the semigroup $S_c(W)$ becomes a flat semiring if the zero element is taken as the top element. If $W$ consists of a single word $w$, then $S_c(W)$ will be denoted by $S_c(w)$.

Another class of flat semirings are hypergraph semirings. Some basic background  and  related concepts on hypergraphs see \cite{Bretto, hamjac, Voloshin}.
Jackson, Zhao et al. in \cite{Jackson} constructed
a flat semiring $S_{\scriptscriptstyle\mathbb{H}}$ from a $k$-uniform hypergraph  $\mathbb{H}$ of girth $m\geq 5$ ($k\geq3$), which is referred as \emph{hypergraph semirings}. The generators of $S_{\scriptscriptstyle\mathbb{H}}$ consist of the multiplicative and additive zero element $0$ and the elements ${\bf a}_v$ for each vertex $v$ of $\mathbb{H}$. The multiplicative part of the semiring is subject to four rules (see, \cite[p. 225]{Jackson}).  The addition on $S_\mathbb{H}$ is defined by the flat semiring addition.
On the other hand,  to study the varieties of semigroup, Jackson et al.in \cite{jackson2} also defined some 3-uniform hypergraph semigroups from 3-uniform hypergraphs of girth at least 4.

Different from \cite{Jackson} and \cite{jackson2}, in this paper we shall introduce  a hypergraph semigroup $SG_{\scriptscriptstyle \mathbb{H}}$ from a 3-hypergraph $\mathbb{H}$ of girth $m\geq 3$
(which need not be either 3-uniform or of girth $m\geq 4$).
We shall prove that each such 3-hypergraph semigroup $SG_{\scriptscriptstyle \mathbb{H}}$ is 0-cancellative. Thus $SG_{\scriptscriptstyle \mathbb{H}}$ becomes a flat semiring by adding the flat semiring addition, which is also written $S_{\scriptscriptstyle \mathbb{H}}$ and called  \emph{3-hypergraph semirings}.
We shall show that all 3-hypergraph semirings $S_{\scriptscriptstyle \mathbb{H}}$ are subdirectly irreducible and nonfinitely based in Section 2.
Also, it will be proved that for a 3-hypergraph $\mathbb{H}$, there exists a 3-uniform hypergraph $\mathbb{H}'$ such that
$\mathbf{V}(S_{\scriptscriptstyle \mathbb{H}})=\mathbf{V}(S_{\scriptscriptstyle{\mathbb{H}'}})$. That is to say, each variety generated by 3-hypergraph semirings is equal to a variety generated by 3-uniform hypergraph semirings.
In Section 3, we shall introduce and study the 2-robustly strong 3-colorable 3-uniform hypergraph, such as 3-uniform $n$-cycle $(n\geq 4)$ and 3-uniform hyperforest,
and hypergraph semirings defined from them.
It will be shown that each variety generated by 2-robustly strong 3-colorable 3-uniform hypergraph semirings is equal to  $\mathbf{V}(S_c(abc))$. Of course, both variety $\mathbf{V}(S_c(abc))$ (see, \cite{Jackson,Ren}) and variety $\mathbf{V}(S_{\scriptscriptstyle \mathbb{H}})$ play key role in the theory of variety of ai-semirings, where 3-uniform hypergraph $\mathbb{H}$ is a 3-cycle.
We shall show that each variety generated by so-called beam-type hypergraph semirings or fan-type hypergraph semirings is equal to the variety $\mathbf{V}(S_{\scriptscriptstyle \mathbb{H}})$ in Section 4. Also, we shall provide an infinite ascending chain in the lattice of subvarieties of the variety generated by all 3-uniform hypergraph semirings.
This implies that the variety generated by all 3-uniform hypergraph semirings has infinitely many subvarieties.

Throughout this paper, 3-hypergraph mentioned is always finite, linear and contains no isolated vertices and loops, and every hyperedge of degree 2
can not be adjacent with other hyperedges.
For other notations and terminologies on semigroup and equational
theory not given in this paper, the reader is referred to
the books \cite{burris}, \cite{Howie} and \cite{Sapir}.

\section{3-hypergraph semigroups and semirings}

In this section we shall introduce and study 3-hypergraph semigroups and 3-hypergraph semirings. Some properties of them are given. For this, the following related concepts on hypergraphs are needed.

Let $V = \{v_1, v_2, \ldots, v_n\}$ be a finite set and $E = \{e_1, e_2, \ldots, e_m\}$ a family of subsets of $V$. The pair $\mathbb{H}=(V, E)$ is called a \emph{hypergraph} with the vertex set $V$ and the hyperedge set $E$ (see, \cite{Bretto, hamjac, Voloshin}).
The elements $v_1, v_2, \ldots, v_n$ are called \emph{vertices} and the sets $e_1, e_2, \ldots, e_m$ are called \emph{hyperedges}. In a hypergraph, two vertices are said to be \emph{adjacent} if there is a hyperedge $e\in E$ that contains both vertices.
Two hyperedges are said to be \emph{adjacent} if their intersection is not empty.
A vertex is called \emph{isolated} if it does not belong to any hyperedge. A hyperedge $e\in E$ such that $|e|=1$ is a \emph{loop}.

Some hyperedges may be the subsets of some other hyperedges; in this case they are called \emph{included}.
A hypergraph is called \emph{simple} if it contains no included hyperedges.
Two simple hypergraph $\mathbb{H}_1$ and $\mathbb{H}_2$ are called \emph{isomorphic} if there exists a one-to-one correspondence between their vertex sets such that any subset of vertices form a  hyperedge in $\mathbb{H}_1$ if and only if the corresponding subset of vertices forms a  hyperedge in $\mathbb{H}_2$.

A subset of a hyperedge is called a \emph{subhyperedge} of $\mathbb{H}$. Two subhyperedges are said to be \emph{linked} if there is a single vertex that completes both of the subhyperedges to full hyperedges of $\mathbb{H}$.
Any hypergraph $\mathbb{H}'= (V', E')$ such that $V'\subseteq V$ and $E'\subseteq E$ is called a  \emph{subhypergraph} of $\mathbb{H}$. In such case, we write $\mathbb{H}'\subseteq \mathbb{H}$.
A hypergraph $\mathbb{H}'= (V', E')$ is called an \emph{induced subhypergraph} of $\mathbb{H}$ if $V'\subseteq V$ and all hyperedges of $\mathbb{H}$ completely contained in $V'$ form the family $E'$. Sometimes we say that $\mathbb{H}'$ is a subhypergraph induced by $V'$.
A hypergraph $\mathbb{H}'= (V', E')$ is called a \emph{partial subhypergraph} of $\mathbb{H}$ if $E'\subseteq E$ and $\cup\{e_i\mid e_i\in E'\}\subseteq V'\subseteq V$. Notice that we may have $V'=V$ or $\cup\{e_i\mid e_i\in E'\}$.

A \emph{cycle} is a sequence $v_0, e_0, v_1, e_1, \ldots, v_{n-1}, e_{n-1}$ alternating between distinct vertices $v_0, \dots, v_{n-1}$ and distinct hyperedges $e_0, \ldots, e_{n-1}$, such that $v_i\in e_i\cap e_{i+1}$ (with addition in the subscript taken modulo $n$), and $n$ is called the length of this cycle. We denote the cycle with the length being equal to $n$ by $n$-cycle.
The \emph{girth} of a hypergraph $\mathbb{H}$ is the length of the shortest cycle if there exists at least a cycle in $\mathbb{H}$, denoted by $g(\mathbb{H})$.
A \emph{hyperforest} means a hypergraph without cycles, which girth is regarded as infinite. A hypergraph is said to be \emph{$n$-cycle hypergraph} if it is exactly a $n$-cycle.

The number $|e_i|$ is called the \emph{degree} of the hyperedge $e_i$.
The rank of a hypergraph $\mathbb{H}$ is $r(\mathbb{H})=max\{|e_i|\mid e_i\in E\}$.
A hypergraph (resp., hyperforest) is called \emph{$r$-hypergraph} (resp., \emph{$r$-hyperforest}) if its rank is $r$.
If $|e|=k$ for all $e\in E$, then $\mathbb{H}$ is a \emph{k-uniform hypergraph}.
A hypergraph is \emph{linear} if any two of its hyperedges intersect in at most one vertex.
Note that for a hypergraph without loops, we always have the following property.
\begin{prop}
Let $\mathbb{H}$ be a hypergraph without loops. Then $\mathbb{H}$ is linear  if and only if  $g(\mathbb{H})\geq 3$. 
\end{prop}
\pf
As $g(\mathbb{H})\geq 2$, we only need to prove that  $\mathbb{H}$ is not linear if and only if $g(\mathbb{H})=2$. 
For necessity, suppose that $\mathbb{H}$ is not linear.  It follows that there exists two distinct hyperedges $e_1$ and $e_2$ in $\mathbb{H}$ such that $|e_1\cap e_2|\geq 2$. Assume that $e_1\cup e_2=\{v_1,\ldots, v_m\}$, $e_1= \{v_1,\ldots, v_k, v_m\}$ and $e_2= \{v_1, v_{k+1}\ldots, v_{m-1}, v_m\}$. It is clear that the sequence $v_1, e_1, v_m, e_2$ is a 2-cycle. This implies that $g(\mathbb{H})=2$, since $\mathbb{H}$ has no loops.

For sufficiency, suppose that $g(\mathbb{H})=2$. Then there exists a 2-cycle $v_1, e_1, v_m, e_2$, where $e_1\neq e_2$ and $v_1\neq v_m$. Thus $\{v_1, v_m\}\subseteq e_1\cap e_2$. This means that $\mathbb{H}$ is not linear. 
\epf

In the following we shall introduce and study  the 3-hypergraph semigroups. A \emph{3-hypergraph semigroup} $SG_{\scriptscriptstyle \mathbb{H}}$ defined from 3-hypergraph $\mathbb{H}=(V, E)$
is expressed as $\langle X\mid R \rangle$,
where the generator set $X$ is equal to $\{{\bf a}_{u} \mid u\in V\} \cup \{0 \}$ and
the relation $R$ is equal to the union of the following subsets of $X^+ \times X^+$:
$$  \begin{array}{ll}
\{({\bf a}_{u}{\bf a}_{v}, {\bf a}_{v}{\bf a}_{u})\mid u,v \,\in V\},  & \{({\bf a}_{u_{1}}\ldots{\bf a}_{u_{k}}, {\bf a})\mid \{u_{1}, \ldots, u_{k}\} \in E, \, {\bf a} \neq0\}, \\
\{({\bf a}_{u}{\bf a}_{v}{\bf a}_{h}, 0) \mid \{u,v,h\} \notin E\},   & \{({\bf a}_{u}{\bf a}_{v}, 0) \mid u=v \,\, or\,\, \{u,v\} \,\, is \,\, not \,\, a \,\,subhyperedge\}, \\
 \{(0{\bf a}_{u}, 0), ({\bf a}_{u}0, 0), (00, 0) \mid u \,\in V\},     & \{({\bf a}_{u_i}{\bf a}_{u_j}, {\bf a}_{u_r}{\bf a}_{u_l})\mid \{u_i, u_j\} \,\, and\,\, \{u_r, u_l\}\,\,are \,\,linked \,\}. \\
\end{array}  $$

It is well-known that the semigroup $SG_{\scriptscriptstyle \mathbb{H}}$ is isomorphic to $X^{+}/R^{\#}$, where $R^{\#}$ denotes the semigroup congruence on free semigroup $X^{+}$ generated by $R$. To prove that the semigroup $X^{+}/R^{\#}$ is 0-cancellative, the following lemma is needed for us.

\begin{lem}\label{lem1}
In the semigroup $X^{+}/R^{\#}$ the following are true:
	\begin{itemize}
		\item [$(i)$]
		$\{u_{i}, u_{j}\}$ is not a subhyperedge of $\mathbb{H}$ if and only if ${\bf a}_{u_{i}}{\bf a}_{u_{j}}/R^{\#} = 0/R^{\#}$.
		\item [$(ii)$]
		$\{u_{i}, u_{j}, u_{l}\}$ is not a hyperedge of $\mathbb{H}$ if and only if ${\bf a}_{u_{i}}{\bf a}_{u_{j}}{\bf a}_{u_{l}}/R^{\#} = 0/R^{\#}$.
		\item [$(iii)$]
		${\bf a}_{u_{1}}{\bf a}_{u_{2}}\cdots{\bf a}_{u_{l}}/R^{\#} = 0/R^{\#}$ $(l\geq 4)$.
		\item [$(iv)$]
		Let $\{u_{1}, \ldots, u_{i}\}$ and $\{v_{1}, \ldots, v_{j}\}$ be distinct subhyperedges $e'_1$ and $e'_2$ of $\mathbb{H}$, respectively. Then ${\bf a}_{u_{1}}\cdots{\bf a}_{u_{i}}/R^{\#} = {\bf a}_{v_{1}}\cdots{\bf a}_{v_{j}}/R^{\#}$ if and only if both $e'_1$ and $e'_2$ are hyperedges, or they are linked.
	\end{itemize}
\end{lem}
\pf
Clearly, $(i)$ and $(ii)$ are valid. \par
$(iii)$ We only need to show that there exists at least one 3-element subset $A$ of $\{u_{1},u_{2}\cdots u_{l}\}$ such that $A$ is not a hyperedge ($l\geq 4$).
Suppose that all 3-element subsets of $\{u_{1},u_{2}\cdots u_{l}\}$ are hyperedges. Thus $\{u_1, u_2, u_3\}$ and $\{u_2, u_3, u_4\}$ are hyperedges, a contradiction.   \par
$(iv)$ Supose ${\bf a}_{u_{1}}\cdots{\bf a}_{u_{i}}/R^{\#} = {\bf a}_{v_{1}}\cdots{\bf a}_{v_{j}}/R^{\#}$.
Since the degree of a hyperedge in $\mathbb{H}$ is 2 or 3, without loss of generality, we have following six cases.

{\bfseries Case 1.}  $e'_1=\{u_{1}\}$ and $e'_2=\{v_{1}\}$. Then ${\bf a}_{u_{1}}/R^{\#}={\bf a}_{v_{1}}/R^{\#}$, which is contrary to ${\bf a}_{u_{1}}$ and ${\bf a}_{v_{1}}$ being different.
Hence the assumption does not hold in this case.

{\bfseries Case 2.}  $e'_1=\{u_{1}\}$ and $e'_2=\{v_{1}, v_{2}\}$. \par
(1) If $e'_2$ is a hyperedge, then ${\bf a}_{v_{1}}{\bf a}_{v_{2}}/R^{\#}={\bf a}/R^{\#}={\bf a}_{u_{1}}/R^{\#}$. Thus $\{u_{1}\}$ is a hyperedge, a contradiction. Hence the assumption does not hold at this point. \par
(2) If $e'_2$ is a proper subhyperedge, then there exists a vertex $v_3\in V$ such that $\{v_{1}, v_{2}, v_{3}\}$ is a hyperedge. Since ${\bf a}_{u_{1}}/R^{\#}={\bf a}_{v_{1}}{\bf a}_{v_{2}}/R^{\#}$, we have that ${\bf a}_{v_{1}}{\bf a}_{v_{2}}{\bf a}_{v_{3}}/R^{\#}={\bf a}/R^{\#}={\bf a}_{u_{1}}{\bf a}_{v_{3}}/R^{\#}$. Thus $\{u_{1}, v_{3}\}$ is a hyperedge, a contradiction. Hence the assumption does not hold  at this point.

{\bfseries Case 3.}  $e'_1=\{u_{1}\}$ and $e'_2=\{v_{1}, v_{2}, v_{3}\}$. Clearly, $e'_2$ is a hyperedge. Then ${\bf a}_{v_{1}}{\bf a}_{v_{2}}{\bf a}_{v_{3}}/R^{\#} = {\bf a}/R^{\#} = {\bf a}_{u_{1}}/R^{\#}$. Thus $\{u_{1}\}$ is a hyperedge, a contradiction. Hence the assumption does not hold in this case.

{\bfseries Case 4.}  $e'_1=\{u_{1}, u_{2}\}$ and $e'_2=\{v_{1}, v_{2}\}$.\par
(1) If $e'_1$ is a hyperedge, then ${\bf a}_{u_{1}}{\bf a}_{u_{2}}/R^{\#}={\bf a}/R^{\#}={\bf a}_{v_{1}}{\bf a}_{v_{2}}/R^{\#}$. Then $e'_2$ is also a hyperedge. \par
(2) If $e'_1$ is a proper subhyperedge, then there exists a vertex $u_3\in V$ such that $\{u_{1}, u_{2}, u_{3}\}$ is a hyperedge. Since ${\bf a}_{u_{1}}{\bf a}_{u_{2}}/R^{\#}={\bf a}_{v_{1}}{\bf a}_{v_{2}}/R^{\#}$, we have that ${\bf a}_{u_{1}}{\bf a}_{u_{2}}{\bf a}_{u_{3}}/R^{\#}={\bf a}/R^{\#}={\bf a}_{v_{1}}{\bf a}_{v_{2}}{\bf a}_{u_{3}}/R^{\#}$. Thus $\{v_{1}, v_{2}, u_{3}\}$ is a hyperedge. Suppose that  $|e'_1\cap e'_2|=1$. Without loss of generality we assume that $u_1=v_1$. Then $\{u_1, u_2, u_3\}$ and $\{u_1, v_2, u_3\}$ are hyperedges, a contradiction.
Hence $e'_1$ and $e'_2$ are linked.

{\bfseries Case 5.}  $e'_1=\{u_{1}, u_{2}\}$ and $e'_2=\{v_{1}, v_{2}, v_{3}\}$. Clearly, $e'_2$ is a hyperedge. Then ${\bf a}_{v_{1}}{\bf a}_{v_{2}}{\bf a}_{v_{3}}/R^{\#}={\bf a}/R^{\#}={\bf a}_{u_{1}}{\bf a}_{u_{2}}/R^{\#}$. Thus $e'_1$ is also a hyperedge.

{\bfseries Case 6.}  $e'_1=\{u_{1}, u_{2}, u_{3}\}$ and $e'_2=\{v_{1}, v_{2}, v_{3}\}$. Clearly, $e'_1$ and $e'_2$ are hyperedges.

The converse half is clear.
\epf

\begin{prop}\label{lem2}
The semigroup $X^{+}/R^{\#}$ is 0-cancellative.
\end{prop}
\pf
Let $a,b,c \in X^{+}$. Suppose that $(a/R^{\#}) (b/R^{\#}) = (a/R^{\#}) (c/R^{\#}) \neq 0/R^{\#}$. Since $a/R^{\#}\neq 0/R^{\#}$ and $(a/R^{\#}) (b/R^{\#})\neq 0/R^{\#}$, it follows from Lemma \ref{lem1} ($i-iii$) that there exists a subhyperedge $\{u_{1}, u_{2}\}$ such that $a/R^{\#}={\bf a}_{u_{1}}{\bf a}_{u_{2}}/R^{\#}$ or there exists a vertex $u_{1}\in V$ such that $a/R^{\#}={\bf a}_{u_{1}}/R^{\#}$.
Similarly, there exists a subhyperedge $\{u_{3}, u_{4}\}$ such that $b/R^{\#}={\bf a}_{u_{3}}{\bf a}_{u_{4}}/R^{\#}$ or there exists a vertex $u_{3}\in V$ such that $b/R^{\#}={\bf a}_{u_{3}}/R^{\#}$;
there exists a subhyperedge $\{u_{5}, u_{6}\}$ such that $c/R^{\#}={\bf a}_{u_{5}}{\bf a}_{u_{6}}/R^{\#}$ or there exists a vertex $u_{5}\in V$ such that $c/R^{\#}={\bf a}_{u_{5}}/R^{\#}$. Without loss of generality, we only consider the following six cases.

{\bfseries Case 1.} $a/R^{\#}={\bf a}_{u_{1}}/R^{\#}$, $b/R^{\#}={\bf a}_{u_{3}}/R^{\#}$ and $c/R^{\#}={\bf a}_{u_{5}}/R^{\#}$. If $u_{3}=u_{5}$, then $b=c$ and so $b/R^{\#} = c/R^{\#}$ as required.
Otherwise, since ${\bf a}_{u_{1}}{\bf a}_{u_{3}}/R^{\#}={\bf a}_{u_{1}}{\bf a}_{u_{5}}/R^{\#}\neq 0/R^{\#}$, it follows from  Lemma \ref{lem1} ($i$) and ($iv$) that subhyperedges $\{u_{1}, u_{3}\}$ and $\{u_{1}, u_{5}\}$ are hyperedges or linked, a contradiction.

{\bfseries Case 2.} $a/R^{\#}={\bf a}_{u_{1}}/R^{\#}$, $b/R^{\#}={\bf a}_{u_{3}}/R^{\#}$ and $c/R^{\#}={\bf a}_{u_{5}}{\bf a}_{u_{6}}/R^{\#}$.
It follows from Lemma \ref{lem1} ($ii$) that $\{u_{1}, u_{5}, u_{6}\}$ is a hyperedge, since ${\bf a}_{u_{1}}{\bf a}_{u_{5}}{\bf a}_{u_{6}}/R^{\#}\neq 0/R^{\#}$. Thus $\{u_{1}, u_{3}\}$ is also a hyperedge, since ${\bf a}_{u_{1}}{\bf a}_{u_{5}}{\bf a}_{u_{6}}/R^{\#}={\bf a}/R^{\#}={\bf a}_{u_{1}}{\bf a}_{u_{3}}/R^{\#}$, a contradiction. Hence the assumption does not hold in this case.

{\bfseries Case 3.} $a/R^{\#}={\bf a}_{u_{1}}/R^{\#}$, $b/R^{\#}={\bf a}_{u_{3}}{\bf a}_{u_{4}}/R^{\#}$ and $c/R^{\#}={\bf a}_{u_{5}}{\bf a}_{u_{6}}/R^{\#}$.
If $\{u_3, u_4\}=\{u_5, u_6\}$, then $b=c$ and so $b/R^{\#} = c/R^{\#}$ as required.
Otherwise, since ${\bf a}_{u_{1}}{\bf a}_{u_{3}}{\bf a}_{u_{4}}/R^{\#}={\bf a}_{u_{1}}{\bf a}_{u_{5}}{\bf a}_{u_{6}}/R^{\#}\neq 0/R^{\#}$, it follows from Lemma \ref{lem1} ($ii$) that $\{u_{1}, u_{3}, u_{4}\}$ and $\{u_{1}, u_{5}, u_{6}\}$ are hyperedges.
Suppose that  $|\{u_{3}, u_{4}\}\cap \{u_{5}, u_{6}\}|=1$. Without loss of generality we assume that $u_3=u_5$. Then $\{u_1, u_3, u_4\}$ and $\{u_1, u_3, u_6\}$ are hyperedges, a contradiction.
Thus $\{u_{3}, u_{4}\}$ and $\{u_{5}, u_{6}\}$ are linked. Hence $b/R^{\#}={\bf a}_{u_{3}}{\bf a}_{u_{4}}/R^{\#}={\bf a}_{u_{5}}{\bf a}_{u_{6}}/R^{\#}=c/R^{\#}$.

{\bfseries Case 4.} $a/R^{\#}={\bf a}_{u_{1}}{\bf a}_{u_{2}}/R^{\#}$, $b/R^{\#}={\bf a}_{u_{3}}/R^{\#}$ and $c/R^{\#}={\bf a}_{u_{5}}/R^{\#}$. If $u_{3}=u_{5}$, then $b=c$ and so $b/R^{\#} = c/R^{\#}$ as required. Otherwise, since ${\bf a}_{u_{1}}{\bf a}_{u_{2}}{\bf a}_{u_{3}}/R^{\#}={\bf a}_{u_{1}}{\bf a}_{u_{2}}{\bf a}_{u_{5}}/R^{\#}\neq 0/R^{\#}$, it follows from Lemma \ref{lem1} ($ii$) that $\{u_{1}, u_{2}, u_{3}\}$ and $\{u_{1}, u_{2}, u_{5}\}$ are hyperedges, a contradiction.

{\bfseries Case 5.}  $a/R^{\#}={\bf a}_{u_{1}}{\bf a}_{u_{2}}/R^{\#}$, $b/R^{\#}={\bf a}_{u_{3}}/R^{\#}$ and $c/R^{\#}={\bf a}_{u_{5}}{\bf a}_{u_{6}}/R^{\#}$. Then the assumption does not hold in this case, since $(a/R^{\#}) (c/R^{\#}) = 0/R^{\#}$.

{\bfseries Case 6.} $a/R^{\#}={\bf a}_{u_{1}}{\bf a}_{u_{2}}/R^{\#}$, $b/R^{\#}={\bf a}_{u_{3}}{\bf a}_{u_{4}}/R^{\#}$ and $c/R^{\#}={\bf a}_{u_{5}}{\bf a}_{u_{6}}/R^{\#}$. Then assumption does not hold in this case, since $(a/R^{\#})(b/R^{\#})=(a/R^{\#})(c/R^{\#})=0/R^{\#}$.

Dually, we can prove that $b/R^{\#} = c/R^{\#}$ if $(b/R^{\#})(a/R^{\#}) = (c/R^{\#})(a/R^{\#}) \neq 0/R^{\#}$.
\epf

Now we can see from \cite[Lemma 2.2]{Jackson} that the semigroup $X^{+}/R^{\#}$ equipped with the flat semiring addition becomes a flat semiring, which is written $S_{\scriptscriptstyle \mathbb{H}}$ and called  \emph{3-hypergraph semiring}.

\begin{prop}\label{prop2.3}
	All 3-hypergraph semirings $S_{\scriptscriptstyle \mathbb{H}}$ are nonfinitely based and subdirectly irreducible.
\end{prop}
\pf
It is easy to see that $S_c(a_1a_2a_3)\in \mathbf{V}(S_{\scriptscriptstyle \mathbb{H}})$, since $S_c(a_1a_2a_3)$ is, up to isomorphism, a subsemiring of $S_{\scriptscriptstyle \mathbb{H}}$.
They are also clear that $S_{\scriptscriptstyle \mathbb{H}}$  is a finite flat semiring with the index 2 and that all noncyclic elements in $S_{\scriptscriptstyle \mathbb{H}}$ form an order ideal of $S_{\scriptscriptstyle \mathbb{H}}$.  Thus
by \cite[Theorem 4.9]{Jackson}, we have directly that $S_{\scriptscriptstyle \mathbb{H}}$ is nonfinitely based.

From \cite[Proposition 2.1]{Wu}, it follows that $S_{\scriptscriptstyle \mathbb{H}}$ is subdirectly irreducible, since $S_{\scriptscriptstyle \mathbb{H}}$ is 2-nil ( i.e., its multiplicative reduct is a $2$-nil semigroup) and has a unique non-zero annihilator $\bar{{\bf a}}$, where $\bar{{\bf a}}$ denotes ${\bf a}/R^{\#}$.
\epf

\begin{lem} \label{lem2.4}
For every 3-hypergraph $\mathbb{H}$, there exists a 3-uniform hypergraph $\mathbb{H}'$ such that
	$\mathbf{V}(S_{\scriptscriptstyle \mathbb{H}})=\mathbf{V}(S_{\scriptscriptstyle{\mathbb{H}'}})$.
\end{lem}
\pf
This statement is true when $\mathbb{H}$ is a 3-uniform hypergraph. Otherwise, suppose that $\mathbb{H}=(V, E)$ with $V=\{u_1, u_2, \ldots, u_n, u_{n+1}, u_{n+2}\}$ is not 3-uniform and $e_{m}=\{u_{n+1}, u_{n+2}\}$ is a hyperedge of $\mathbb{H}$ with degree 2.
Write  $\mathbb{H}_1$ to denote the subhypergraph of $\mathbb{H}$ induced by the vertex set $V\setminus e_{m}$.
In the following we shall show that $\mathbf{V}(S_{\scriptscriptstyle \mathbb{H}}) = \mathbf{V}(S_{\scriptscriptstyle\mathbb{H}_1})$.
It is easy to see that  $S_{\scriptscriptstyle \mathbb{H}_1}\in \mathbf{V}(S_{\scriptscriptstyle \mathbb{H}})$ and so $\mathbf{V}(S_{\scriptscriptstyle \mathbb{H}_1}) \subseteq \mathbf{V}(S_{\scriptscriptstyle \mathbb{H}})$.

On the other hand, it follows that $\mathbb{H}_1$ must have a hyperedge of degree 3, say  $\{u_1, u_2, u_3\}$ without loss of generality,
since every hyperedge of degree 2 in $\mathbb{H}$ can not be adjacent with a hyperedge of degree 2 or 3. Thus 3-hypergraph semiring $S_{\scriptscriptstyle \mathbb{H}_1}$ is generated by $\{\bar{{\bf a}}_{u_1}, \bar{{\bf a}}_{u_2}, \ldots, \bar{{\bf a}}_{u_n} \}$, where $\bar{{\bf a}}_{u_i}$ denotes ${\bf a}_{u_i}/R^{\#}$.

Let $(S_{{\scriptscriptstyle \mathbb{H}}_1})^2$ denote the direct product of 2 copies of $S_{{\scriptscriptstyle \mathbb{H}}_1}$. Consider the subsemiring $A$ of $(S_{{\scriptscriptstyle \mathbb{H}}_1})^2$ generated by $\{\alpha_1, \alpha_2, \ldots, \alpha_n, , \alpha_{n+1}, \alpha_{n+2}\}$, where
$$\alpha_1=(\bar{{\bf a}}_{u_1}, \bar{{\bf a}}_{u_1}), \alpha_2=(\bar{{\bf a}}_{u_2}, \bar{{\bf a}}_{u_2}),
\alpha_3=(\bar{{\bf a}}_{u_3}, \bar{{\bf a}}_{u_3}),
\ldots, \alpha_n=(\bar{{\bf a}}_{u_n}, \bar{{\bf a}}_{u_n}), $$
$$\alpha_{n+1}=(\bar{{\bf a}}_{u_1}\bar{{\bf a}}_{u_2}, \bar{{\bf a}}_{u_3}), \alpha_{n+2}=(\bar{{\bf a}}_{u_3}, \bar{{\bf a}}_{u_1}\bar{{\bf a}}_{u_2}).$$
Take $J=\{(x_1, x_2)\in A\mid (\exists t\in \{1, 2\})~ x_t=0 \}$.
It is easy to verify that $J$ is an ideal of $A$ and $A/J$ is isomorphic to $S_{\scriptscriptstyle \mathbb{H}}$ under the obvious map defined on generators by $\bar{\alpha}_i\mapsto \bar{{\bf a}}_{u_i}$, where $\bar{\alpha}_i$ denote $\alpha_i/J$. Thus $S_{\scriptscriptstyle \mathbb{H}}\in \mathbf{V}(S_{{\scriptscriptstyle \mathbb{H}}_1})$ and so $\mathbf{V}(S_{\scriptscriptstyle \mathbb{H}}) \subseteq \mathbf{V}(S_{{\scriptscriptstyle \mathbb{H}}_1})$. Now we have  shown that  $\mathbf{V}(S_{\scriptscriptstyle \mathbb{H}}) = \mathbf{V}(S_{{\scriptscriptstyle \mathbb{H}}_1})$. We may set $\mathbb{H}' = \mathbb{H}_1$
if $\mathbb{H}_1$ is 3-uniform. This completes our proof.
Otherwise, suppose that $\mathbb{H}_1=(V_1, E_1)$ is not 3-uniform and $e_{m-1}$ is a hyperedge of $\mathbb{H}_1$ with degree 2. Write  $\mathbb{H}_2$ to denote the subhypergraph of $\mathbb{H}_1$
induced by the vertex set $V_1\setminus e_{m-1}$.
Similarly, we can get that $\mathbf{V}(S_{{\scriptscriptstyle \mathbb{H}}_1}) = \mathbf{V}(S_{{\scriptscriptstyle \mathbb{H}}_2})$ and so  $\mathbf{V}(S_{\scriptscriptstyle \mathbb{H}}) = \mathbf{V}(S_{{\scriptscriptstyle \mathbb{H}}_2})$. We may set $\mathbb{H}' = \mathbb{H}_2$
if $\mathbb{H}_2$ is 3-uniform. This completes our proof. Otherwise, we are going on the above procedure. At most steps, we can see that there exists a 3-uniform hypergraph $\mathbb{H}_n = (V_n, E_n)$
(i.e., the subhypergraph of $\mathbb{H}$ induced by the vertex set $V_n = \cup\{e \mid (\exists~e\in E)~ |e|=3 \}$)
such that $\mathbf{V}(S_{\scriptscriptstyle \mathbb{H}})=\mathbf{V}(S_{{\scriptscriptstyle \mathbb{H}}_{n}})$. Thus we may set $\mathbb{H}' = \mathbb{H}_n$. This completes our proof.
\epf

\begin{thm}\label{thm2.5} Each variety generated by 3-hypergraph semirings is equal to a variety generated by 3-uniform hypergraph semirings.
\end{thm}
\pf
By Lemma \ref{lem2.4} we have that for any given the family $\{\mathbb{H}_i \mid i \in I\}$ of 3-hypergraphs,
there exists a family $\{\mathbb{H}^\prime_i \mid i \in I\}$ of 3-uniform hypergraphs such that  $ \mathbf{V}(S_{{\scriptscriptstyle \mathbb{H}}_i}) = \mathbf{V}(S_{{\scriptscriptstyle \mathbb{H}}^\prime_i})$ for all $i \in I$.
This implies that
\begin{equation*}
\mathbf{V}({\displaystyle\bigcup}_{i \in I} S_{{\scriptscriptstyle \mathbb{H}}_i})
= {\displaystyle \bigvee}_{i \in I} \mathbf{V}(S_{{\scriptscriptstyle \mathbb{H}}_i})
={\displaystyle \bigvee}_{i \in I} \mathbf{V}(S_{{\scriptscriptstyle \mathbb{H}}^\prime_i})
=\mathbf{V}({\displaystyle\bigcup}_{i \in I} S_{{\scriptscriptstyle \mathbb{H}}^\prime_i}).  \qedhere
\end{equation*}
\epf
The above theorem tells us that
we only need to consider the varieties generated by 3-uniform hypergraph semirings, to study the varieties generated by  3-hypergraph semirings.
Thus the 3-hypergraphs mentioned in remainder of the paper are always 3-uniform.

\section{Varieties generated by 2-robustly strong 3-colorable 3-hypergraph semirings}

In this section we shall introduce and study the 2-robustly strong 3-colorable 3-hypergraph and 3-hypergraph semirings defined from such hypergraph. Some examples of 2-robustly strong 3-colorable 3-hypergraphs are given.
We also show that each variety generated by such hypergraph semirings is equal to  $\mathbf{V}(S_c(abc))$.

Recall that a \emph{$3$-coloring} of a hypergraph $\mathbb{H} = (V, E)$ is a mapping $\varphi$ from the vertex set $V$ to the color set $\{0, 1, 2\}$ such that every vertex has just one color and no hyperedge with a cardinality more than 1 is monochromatic.
A 3-coloring of $\mathbb{H}$ is said to be \emph{strong } if every hyperedge of $\mathbb{H}$ is polychromatic (i.e., has all vertices colored differently).
Also, a hypergraph $\mathbb{H}$ is said to be \emph{strong 3-colorable} if it has a strong 3-coloring (see, \cite{Bretto}).

We say that a hypergraph $\mathbb{H}$ is \emph{2-robustly strong 3-colorable} if for any given $2$-element subset $\{u, v\}$ of $V$, every valid partial assignment $\varphi'$ from $\{u, v\}$ to $\{0, 1, 2\}$ (i.e., $\varphi'(u)\neq \varphi'(v)$ when $\{u, v\}$ is a subhyperedge of $\mathbb{H}$) can extend to a full strong 3-coloring.
It is clear that a 2-robustly strong 3-colorable hypergraph is strong 3-colorable. The following example tells us that not every strong 3-colorable hypergraph is 2-robust.

\begin{exa}\label{exa3.1}
A 3-cycle 3-hypergraph $\mathbb{H}$ is strong 3-colorable but not 2-robustly strong 3-colorable.
\end{exa}
\pf
Assume that $\mathbb{H}=(V, E)$ is the sequence $u_1, e_1, u_3, e_2, u_5, e_3$, where
$$e_1 = \{u_1, u_2, u_3\}, e_2 = \{u_3, u_4, u_5\}, e_3 = \{u_5, u_6, u_1\}.$$
Define $\varphi: V\rightarrow \{0, 1, 2\}$ by
$$\varphi(u_1)=\varphi(u_4)=1,~\varphi(u_2)=\varphi(u_5)=2,~\varphi(u_3)=\varphi(u_6)=0.$$
Then $\varphi$ is a strong 3-coloring of $\mathbb{H}$ and so $\mathbb{H}$ is strong 3-colorable.

Given a valid partial assignment $\varphi'$ from $\{u_1, u_4\}$ to $\{0, 1, 2\}$ by $u_1 \mapsto 0, u_4 \mapsto 1$, we assume that there exists  a strong 3-coloring $\varphi$ of $\mathbb{H}$ which  extends $\varphi'$. Then $\varphi(u_3)=2$ and so $\varphi(u_5)=0$, a contradiction. Thus $\mathbb{H}$ is not $2$-robustly strong 3-colorable.
\epf

Recall that a \emph{leaf} in a hypergraph is a hyperedge that is has just at most one vertex in common with any other hyperedge. It is well-konwn that every finite hyperforest either contains no hyperedges at all, or contains a leaf (see, \cite{Jackson}).
Also, it is easy to see that for a strong 3-colorable hypergraph $\mathbb{H}$, every partial assignament of any one vertex can extend to a full strong 3-coloring of $\mathbb{H}$.
In fact, we have

\begin{exa}
A 3-hyperforest $\mathbb{H}$ is 2-robustly strong 3-colorable.
\end{exa}
\pf
Let $\varphi'$ be an arbitrary valid partial assignment from $\{u, v\}$ to $\{0, 1, 2\}$ for any $2$-element subset $\{u, v\}$ of the vertex set $V$.
We shall prove by induction on the number of hyperedges of $\mathbb{H}=(V, E)$, denoted by $m(\mathbb{H})$, that  $\varphi'$  can extend to a strong 3-coloring of $\mathbb{H}$.

If $m(\mathbb{H})=1$, then the statement is obviously true.
We assume that the statement is true when $m(\mathbb{H})=n$, and consider the case of $m(\mathbb{H})=n+1$.
Let $e$ be a leaf of $\mathbb{H}$. Then  $e \cap e'= \emptyset$ for each $e'\in E \setminus \{e\}$ or there exists a vertex $w$
such that  $e\cap e'=\{w\}$ for each $e'\in E \setminus \{e\}$ with $e \cap e' \neq \emptyset$.

(1) Suppose that $e \cap e'= \emptyset$ for each $e'\in E \setminus \{e\}$. Let  $\mathbb{H}'=(V', E')$  denote the subhypergraph of $\mathbb{H}$ induced
by the vertex set $V\setminus e$.
It is easy to see that if $u, v\in e$, then $\varphi'$ can extend to a strong 3-coloring of $\mathbb{H}$,
since $\mathbb{H}'$ is strong 3-colorable.
Otherwise, $u, v\in V'$, or $u\in V'$ and $ v\in e$.
Assume that $u, v\in V'$ ($u\in V'$ and $ v\in e$, respectively). Then there exists a strong 3-coloring $\varphi^*$ of $\mathbb{H}'$ such that $\varphi^* (u) = \varphi' (u)$ and
$\varphi^* (v) = \varphi' (v)$ ($\varphi^* (u) = \varphi' (u)$, respectively), since it is clear by the induction hypothesis that
$\mathbb{H}'$ is 2-robustly strong 3-colorable. This implies that $\varphi'$ can extend to a strong 3-coloring of $\mathbb{H}$.

(2) Suppose that $e\cap e'=\{w\}$ for some $e'\in E \setminus \{e\}$.
Let $\mathbb{H}'=(V', E')$ denote the partial subhypergraph of $\mathbb{H}$ where $E'=E \setminus \{e\}$ and $V'= \cup\{e_i\mid e_i\in E'\}$.
It is easy to see that if $u, v\in e \setminus \{w\}$, then $\varphi'$ can extend to a strong 3-coloring of $\mathbb{H}$, since $\mathbb{H}'$ is strong 3-colorable.
Otherwise, $u, v\in V'$, or $u\in V'$ and $ v\in e\setminus \{w\}$.
Assume that $u, v\in V'$ ($u\in V'$ and $ v\in e\setminus \{w\}$, respectively). Then there exists a strong 3-coloring $\varphi^*$ of $\mathbb{H}'$ such that $\varphi^* (u) = \varphi' (u)$ and
$\varphi^* (v) = \varphi' (v)$ ($\varphi^* (u) = \varphi' (u)$, respectively), since it is well known by the induction hypothesis that
$\mathbb{H}'$ is 2-robustly strong 3-colorable. This implies that $\varphi'$ can extend to a strong 3-coloring of $\mathbb{H}$.
\epf
Another example of 2-robustly strong 3-colorable hypergraph is as follows:
\begin{exa}
A $n$-cycle 3-hypergraph $\mathbb{H}_n$ $(n\geq 4)$ is 2-robustly strong 3-colorable.
\end{exa}
\pf
(1)	When $n=4$, we may assume that 4-cycle $\mathbb{H}_4=(V_4, E_4)$ is a sequence $u_1, e_1, u_3, e_2, u_5, e_3, u_7, e_4,$ where
$$e_1 = \{u_1, u_2, u_3\}, e_2 = \{u_3, u_4, u_5\}, e_3 = \{u_5, u_6, u_7\}, e_4 = \{u_7, u_8, u_1\}.$$

It will be needed to prove that for any given $2$-element subset $\{u, v\}$ of $V_4$, every valid partial assignment $\varphi'$ from $\{u, v\}$ to $\{0, 1, 2\}$  can extend to a full strong 3-coloring.

(i) If $\{u, v\}$ is a subhyperedge, then we may assume without loss of generality that $\{u, v\}\subseteq e_1$.
Note that the valid partial assignment $\varphi'$ from $\{u_1, u_2\}$ to $\{0, 1, 2\}$ defined by $\varphi' (u_1) = \varphi'' (u_1)$ and $\varphi' (u_2) = \{0, 1, 2\} \setminus \{\varphi'' (u_1), \varphi'' (u_3)\}$ and a valid partial assignment $\varphi''$ from $\{u_1, u_3\}$ to $\{0, 1, 2\}$
defined by $\varphi'' (u_1) = \varphi' (u_1)$ and $\varphi'' (u_3) = \{0, 1, 2\} \setminus \{\varphi' (u_1), \varphi' (u_2)\}$
are determined by each other.
So we only need to consider the case of $\{u, v\} = \{u_1, u_2\}$ or $\{u_2, u_3\}$.

If $\{u, v\} = \{u_1, u_2\}$, then we suppose without loss of generality that $\varphi' (u_1) =0 $ and $\varphi' (u_2) = 1$.
Thus the full strong 3-coloring $\varphi_1: V_4\rightarrow \{0, 1, 2\}$ defined by
$$\varphi_1(u_1)=\varphi_1(u_5)=0,~\varphi_1(u_2)=\varphi_1(u_4)=\varphi_1(u_6)=\varphi_1(u_8)=1,~\varphi_1(u_3)=\varphi_1(u_7)=2$$
is a extension of $\varphi'$.

Dually, we can prove that every valid partial assignment $\varphi'$ from $\{u_2, u_3\}$ to $\{0, 1, 2\}$  can extend to a full strong 3-coloring.
Therefore,  If $\{u, v\}$ is a subhyperedge, then  every valid partial assignment $\varphi'$ from $\{u, v\}$ to $\{0, 1, 2\}$  can extend to a full strong 3-coloring.

(ii) If $\{u, v\}$ is not a subhyperedge, then for $u=u_1$, we only need to consider that $\{u, v\} =\{u_1, u_4\}, \{u_1, u_6\}, \{u_1, u_5\}$, respectively.

{\bfseries Case 1.} $\{u, v\} =\{u_1, u_4\}$.
Assume without loss of generality that $\varphi' (u_1) = 0$ and $ \varphi' (u_4) = 1$, or  $\varphi' (u_1) = 0$ and $ \varphi' (u_4) = 0$.
If $\varphi' (u_1) = 0$ and $ \varphi' (u_4) = 1$, then the full strong 3-coloring $\varphi_1$ is a extension of $\varphi'$.
If $\varphi' (u_1) = 0$ and $ \varphi' (u_4) = 0$, then the full strong 3-coloring $\varphi_2: V_4\rightarrow \{0, 1, 2\}$ by
$$\varphi_2(u_1)=\varphi_2(u_4)=\varphi_2(u_6)=0,~\varphi_2(u_2)=\varphi_2(u_5)=\varphi_2(u_8)=1,~\varphi_2(u_3)=\varphi_2(u_7)=2 $$
is a extension of $\varphi'$.

Dually, we can show that every valid partial assignment $\varphi'$ from $\{u_1, u_6\}$ to $\{0, 1, 2\}$  can extend to a full strong 3-coloring.

{\bfseries Case 2.} $\{u, v\} =\{u_1, u_5\}$.
Assume without loss of generality that $\varphi' (u_1) = 0$ and $ \varphi' (u_5) = 1$, or  $\varphi' (u_1) = 0$ and $ \varphi' (u_5) = 0$.
If $\varphi' (u_1) = 0$ and $ \varphi' (u_5) = 1$, then the full strong 3-coloring $\varphi_2$ is a extension of $\varphi'$.
If $\varphi' (u_1) = 0$ and $ \varphi' (u_5) = 0$, then the full strong 3-coloring $\varphi_1$ is a extension of $\varphi'$.

Thus if $u=u_1$, then every valid partial assignment $\varphi'$ from $\{u, v\}$ to $\{0, 1, 2\}$  can extend to a full strong 3-coloring.
Dually, if $u\in \{u_3, u_5, u_7\}$, then  $\varphi'$ can extend to a full strong 3-coloring.

If $u=u_2$, then we only need to consider that $\{u, v\} =\{u_2, u_4\}, \{u_2, u_8\}, \{u_2, u_6\}$, respectively.

{\bfseries Case 1.} $\{u, v\} =\{u_2, u_4\}$.
Assume without loss of generality that $\varphi' (u_2) = 0$ and $ \varphi' (u_4) = 1$, or  $\varphi' (u_2) = 0$ and $ \varphi' (u_4) = 0$.
If $\varphi' (u_2) = 0$ and $ \varphi' (u_4) = 1$, then the full strong 3-coloring $\varphi_3: V_4\rightarrow \{0, 1, 2\}$ by
$$\varphi_3(u_1)=\varphi_3(u_4)=\varphi_3(u_6)=1,~\varphi_3(u_2)=\varphi_3(u_5)=\varphi_3(u_8)=0,~\varphi_3(u_3)=\varphi_3(u_7)=2 $$
is a extension of $\varphi'$.
If $\varphi' (u_2) = 0$ and $ \varphi' (u_4) = 0$, then the full strong 3-coloring $\varphi_4: V_4\rightarrow \{0, 1, 2\}$ by
$$\varphi_4(u_1)=\varphi_4(u_5)=2,~\varphi_4(u_2)=\varphi_4(u_4)=\varphi_4(u_6)=\varphi_4(u_8)=0,~\varphi_4(u_3)=\varphi_4(u_7)=1 $$
is a extension of $\varphi'$.

Dually, we can show that every valid partial assignment $\varphi'$ from $\{u_2, u_8\}$ to $\{0, 1, 2\}$  can extend to a full strong 3-coloring.

{\bfseries Case 2.} $\{u, v\} =\{u_2, u_6\}$.
Assume without loss of generality that $\varphi' (u_2) = 0$ and $ \varphi' (u_6) = 1$, or  $\varphi' (u_2) = 0$ and $ \varphi' (u_6) = 0$.
If $\varphi' (u_2) = 0$ and $ \varphi' (u_6) = 1$, then the full strong 3-coloring $\varphi_3$ is a extension of $\varphi'$.
If $\varphi' (u_2) = 0$ and $ \varphi' (u_6) = 0$, then the full strong 3-coloring $\varphi_4$ is a extension of $\varphi'$.

Thus if $u=u_2$, then every valid partial assignment $\varphi'$ from $\{u, v\}$ to $\{0, 1, 2\}$  can extend to a full strong 3-coloring.
Dually, if $u\in \{u_4, u_6, u_8\}$, then  $\varphi'$ can extend to a full strong 3-coloring. Therefore,  If $\{u, v\}$ is not a subhyperedge, then  $\varphi'$ can extend to a full strong 3-coloring.

Summarizing above cases we show that $\mathbb{H}_4$ is 2-robustly strong 3-colorable.

(2) Suppose that $\mathbb{H}_n=(V_n, E_n)$ is 2-robustly strong 3-colorable,
where $\mathbb{H}_n$ is a sequence $u_1, e_1, u_3, e_2, \ldots , u_{2n}, e_n,$ and
$$e_1 = \{u_1, u_2, u_3\}, e_2 = \{u_3, u_4, u_5\}, e_3 = \{u_5, u_6, u_7\}, \ldots ,  e_n = \{u_{2n-1}, u_{2n}, u_1\}.$$

We shall prove that $\mathbb{H}_{n+1}=(V_{n+1}, E_{n+1})$ is 2-robustly strong 3-colorable, where $\mathbb{H}_{n+1}$ is a sequence
$u_1, e_1, u_3, e_2, \ldots , u_{2n+2}, e_{n+1},$ and
$$e_1 = \{u_1, u_2, u_3\}, e_2 = \{u_3, u_4, u_5\}, e_3 = \{u_5, u_6, u_7\}, \ldots ,  e_{n+1} = \{u_{2n+1}, u_{2n+2}, u_1\}.$$

It will be needed to prove that for any given $2$-element subset $\{u, v\}$ of $V_{n+1}$, every valid partial assignment from $\{u, v\}$ to $\{0, 1, 2\}$ can extend to a full strong 3-coloring.
Since $\mathbb{H}_{n+1}$ is symmetric, without loss of generality we only need to prove that for any $\{u, v\} \subseteq V_{n+1}\setminus \{u_{2n+2}, u_{2n+1}, u_{2n}, u_{2n-1}\}$, every valid partial assignment $\varphi'$ from $\{u, v\}$ to $\{0, 1, 2\}$ can extend to a strong 3-coloring of $\mathbb{H}_{n+1}$.
Since $\mathbb{H}_n$ is 2-robustly strong 3-colorable, there exists a strong 3-coloring $\varphi_n$  of $\mathbb{H}_{n}$ which extends $\varphi'$. Then the strong 3-coloring $\varphi_{n+1}: V_{n+1}\rightarrow \{0, 1, 2\}$ of $\mathbb{H}_{n+1}$ defined by
\begin{equation*}
\varphi_{n+1}(u_i)= \begin{cases}
\varphi_{n}(u_i) & \text{if} \, u_i\in V_{n+1}\setminus \{u_{2n+2}, u_{2n+1}, u_{2n}\};\\
\varphi_{n}(u_1)  & \text{if} \, u_i = u_{2n};\\
\varphi_{n}(u_{2n})  & \text{if} \, u_i = u_{2n+1};\\
\varphi_{n}(u_{2n-1}) & \text{if} \, u_i = u_{2n+2},
\end{cases}
\end{equation*}
is a extension of $\varphi'$, as required.
\epf

From Example \ref{exa3.1}, we can see that any hypergraph $\mathbb{H}$ with $g(\mathbb{H})= 3$ is not 2-robustly strong 3-colorable. Thus for a 2-robustly strong 3-colorable hypergraph, its girth is at least 4, since the girth of any hypergraph in this paper is always at least 3. Furthermore, we have the following consequence.
\begin{thm}\label{zhuose}
If 3-hypergraph $\mathbb{H} = (V, E)$ is 2-robustly strong 3-colorable, then $\mathbf{V}(S_{\scriptscriptstyle \mathbb{H}}) =  \mathbf{V}(S_c(abc))$.
\end{thm}
\pf
It is easy to see that  $\mathbf{V}(S_{\scriptscriptstyle \mathbb{H}}) \supseteq  \mathbf{V}(S_c(abc))$, since  $S_c(abc)$ is a subsemiring of each 3-hypergraph semiring up to isomorphism. To show the converse inclusion, we need only to prove that $S_{\scriptscriptstyle \mathbb{H}} \in \mathbf{V}(S_c(abc))$.

Let $E'$ denote the set of all subhyperedges of $\mathbb{H}$.
Clearly, $\mathbb{H}$ is  strong 3-colorable and the set $T$ of all strong 3-colorings of $\mathbb{H}$ is finite. Take $A$ to be the subsemiring of $(S_c(abc))^T$ generated by
$$\{\alpha_u\in (S_c(abc))^T \mid u\in V \},$$
where for each $\varphi\in T$,
\begin{equation*}
\alpha_u(\varphi)= \begin{cases}
a & \text{if} \, \varphi(u)=0;\\
b  & \text{if} \, \varphi(u)=1;\\
c & \text{if} \, \varphi(u)=2.
\end{cases}
\end{equation*}
It is clear that $J=\{\alpha\in A\mid (\exists\,  \varphi\in T) \, \alpha(\varphi)=0\}$ is an ideal of the semiring $A$. Also, the quotient semiring $A/J$ is flat. In fact, suppose that $\alpha/J\neq \beta/J$ in $A/J$. Then $\alpha\neq \beta$ and so there exists $\varphi\in T$ such that $\alpha(\varphi)\neq \beta(\varphi)$ in $S_c(abc)$. Thus,  $(\alpha+\beta)(\varphi)=\alpha(\varphi)+\beta(\varphi)=0$, since $S_c(abc)$ is a flat semiring. This implies that  $\alpha+\beta\in J$ and so  $\alpha/J+\beta/J=(\alpha+\beta)/J=0/J$ in $A/J$. This  shows that $A/J$ is flat, as required.

It is also clear that the flat semiring $S_{\scriptscriptstyle \mathbb{H}}$ and $A/J$ are isomorphic if and only if the multiplicative reducts of them are isomorphic.
To show that the multiplicative reducts of flat semirings $S_{\scriptscriptstyle \mathbb{H}}$ and $A/J$ are isomorphic, we need only prove that the multiplicative reduct of $A/J$ is isomorphic to
the semigroup $X^{+}/R^{\#}$,
since  the multiplicative reduct of $S_{\scriptscriptstyle \mathbb{H}}$ is isomorphic to $X^{+}/R^{\#}$.

It is enough to prove that  the extension $\psi$ of the mapping defined on generators by  ${\bf a}_{u}/R^{\#}\mapsto \alpha_{u}/J$ is an isomorphism from $X^{+}/R^{\#}$ to the multiplicative reduct of $A/J$.

For this, the following claim will be needed for us.

\begin{Claim}\label{claim1}
	The following are true in the multiplicative reduct of $A/J$.
	\begin{itemize}
		\item [$(i)$]
		$\{u_{i}, u_{j}\}$ is not a subhyperedge of $\mathbb{H}$ if and only if $\alpha_{u_{i}}\alpha_{u_{j}}/J = 0/J$.
		\item [$(ii)$]
		$\{u_{i}, u_{j}, u_{l}\}$ is not a hyperedge of $\mathbb{H}$ if and only if $\alpha_{u_{i}}\alpha_{u_{j}}\alpha_{u_{l}}/J = 0/J$.
		\item [$(iii)$]
		$\alpha_{u_{1}}\alpha_{u_{2}}\cdots\alpha_{u_{l}}/J = 0/J$ $(l\geq 4)$.
		\item [$(iv)$]
		Let $\{u_{1}, \ldots, u_{i}\}$ and $\{v_{1}, \ldots, v_{j}\}$ be distinct subhyperedges $e'_1$ and $e'_2$ of $\mathbb{H}$, respectively. Then $\alpha_{u_{1}}\cdots\alpha_{u_{i}}/J = \alpha_{v_{1}}\cdots\alpha_{v_{j}}/J$ if and only if $e'_1$ and $e'_2$ are hyperedges, or $e'_1$ and $e'_2$ are linked.
	\end{itemize}
\end{Claim}

\pf
We only need to show that the multiplicative reduct of $A/J$ satisfies the following five rules (1--5), and the rest proof is similar as that in Lemma \ref{lem1}.

(1) $\alpha_u\alpha_v/ J=0/ J$ in the multiplicative reduct of $A/J$ if $u=v$ or $\{u, v\}\notin E'$.
In fact, if $u=v$ or $\{u, v\}\notin E'$, then by 2-robustly strong 3-colorability, there exists a strong 3-coloring $\varphi$ such that $\varphi(u)=\varphi(v)=0$. Thus it follows immediately that
$$\alpha_u\alpha_v(\varphi)=\alpha_u(\varphi)\alpha_v(\varphi)=aa=0,$$
and so $\alpha_u\alpha_v\in J$. This shows that $\alpha_u\alpha_v/ J=0/ J$, as required.

(2) $(\alpha_u/J)(\alpha_v/ J)=(\alpha_v/J)(\alpha_u/ J)$.
Indeed, it is clear that
$$(\alpha_u/J)(\alpha_v/ J)=\alpha_u\alpha_v/ J=\alpha_v\alpha_u/J=(\alpha_v/J)(\alpha_u/ J),$$ since $S_c(abc)$ is commutative.

(3) $\alpha_{u_{i}}\alpha_{u_{j}}\alpha_{u_{l}}/J= 0/J$ if $\{u_{i}, u_{j}, u_{l}\} \notin E$. Indeed, by \cite[Lemma 3.2 (III)]{Jackson}, there exists a 2-element subset of $\{u_{i}, u_{j}, u_{l}\}$, without loss of generality denoted by $\{u_{i}, u_{j}\}$, such that $\{u_{i}, u_{j}\}\notin E'$ since $g(\mathbb{H})\geq 4$. Thus $\alpha_{u_{i}}\alpha_{u_{j}}= 0/J$ and so $\alpha_{u_{i}}\alpha_{u_{j}}\alpha_{u_{l}}/J= 0/J$.

(4) $\alpha_{u_{i}}\alpha_{u_{j}}\alpha_{u_{l}}/J=\alpha_{v_{i}}\alpha_{v_{j}}\alpha_{v_{l}}/J$ for any $\{u_{i}, u_{j}, u_{l}\}, \{v_{i}, v_{j}, v_{l}\} \in E$.
In fact, we have
$$\alpha_{u_{i}}\alpha_{u_{j}}\alpha_{u_{l}}(\varphi) = \alpha_{u_{i}}(\varphi)\alpha_{u_{j}}(\varphi)\alpha_{u_{l}}(\varphi) =abc$$
for each $\varphi\in T$. Dually, we can see that $ \alpha_{v_{i}}\alpha_{v_{j}}\alpha_{v_{l}}(\varphi)=abc$.
That is to say, $\alpha_{u_{i}}\alpha_{u_{j}}\alpha_{u_{l}}=\alpha_{v_{i}}\alpha_{v_{j}}\alpha_{v_{l}}$  and so $\alpha_{u_{i}}\alpha_{u_{j}}\alpha_{u_{l}}/J=\alpha_{v_{i}}\alpha_{v_{j}}\alpha_{v_{l}}/J$.

(5) $\alpha_{u_{i}}\alpha_{u_{j}}/J=\alpha_{v_{i}}\alpha_{v_{j}}/J$ if $ \{u_{i}, u_{j}\}$ and $\{v_{i}, v_{j}\}$ are linked.
Indeed, if $ \{u_{i},u_{j}\}$ and $\{v_{i}, v_{j}\}$ are linked, then there exists $w\in V$ such that $\{u_{i}, u_{j}, w\}, \{v_{i}, v_{j}, w\}\in E$. Thus $\alpha_{u_{i}}\alpha_{u_{j}}\alpha_w/J=\alpha_{v_{i}}\alpha_{v_{j}}\alpha_w/J$ by (4). This implies that $\alpha_{u_{i}}\alpha_{u_{j}}/J=\alpha_{v_{i}}\alpha_{v_{j}}/J$, since the multiplicative reduct of $A/J$ satisfies $0$-cancellative laws.
\epf
Now we may complete the main proof of this proposition.
To show that $\psi$ is well-defined and a monomorphism, we shall prove that
$$\alpha_{u_{1}}\cdots \alpha_{u_{i}}/J=\alpha_{v_{1}}\cdots \alpha_{v_{j}}/J \Longleftrightarrow {\bf a}_{u_{1}}\cdots{\bf a}_{u_{i}}/R^{\#}= {\bf a}_{v_{1}}\cdots{\bf a}_{v_{j}}/R^{\#},$$
for any $\alpha_{u_{1}} \cdots \alpha_{u_{i}}/J, \alpha_{v_{1}} \cdots  \alpha_{v_{j}}/J\in A/J$ and
${\bf a}_{u_{1}} \cdots {\bf a}_{u_{i}}/R^{\#},  {\bf a}_{v_{1}}  \cdots {\bf a}_{v_{j}}/R^{\#} \in X^{+}/R^{\#}$,
where $u_{1}, \ldots, u_{i}, v_{1}, \ldots, v_{j}\in V$.

Firstly, it is clear that if $\{u_{1}, \ldots, u_{i}\}= \{v_{1}, \ldots, v_{j}\}$, then
$$\alpha_{u_{1}}\cdots \alpha_{u_{i}}/J=\alpha_{v_{1}}\cdots \alpha_{v_{j}}/J \Longleftrightarrow {\bf a}_{u_{1}}\cdots{\bf a}_{u_{i}}/R^{\#}= {\bf a}_{v_{1}}\cdots{\bf a}_{v_{j}}/R^{\#}.$$

We now need only to consider the case $\{u_{1}, \ldots, u_{i}\}\neq \{v_{1}, \ldots, v_{j}\}$.
Suppose that $\alpha_{u_{1}}\cdots \alpha_{u_{i}}/J=\alpha_{v_{1}}\cdots \alpha_{v_{j}}/J$. If $\alpha_{u_{1}}\cdots \alpha_{u_{i}}/J =0/J$, then $\alpha_{v_{1}}\cdots \alpha_{v_{j}}/J =0/J$.
Also, it is easy to see from Claim \ref{claim1}  that $\{u_{1}, \ldots, u_{i}\}, \{v_{1}, \ldots, v_{j}\} \notin E'$. Thus by Lemma \ref{lem1} that  ${\bf a}_{u_{1}}\cdots{\bf a}_{u_{i}}/R^{\#}=0/R^{\#}= {\bf a}_{v_{1}}\cdots{\bf a}_{v_{j}}/R^{\#}$, as required. Otherwise, $\alpha_{u_{1}}\cdots \alpha_{u_{i}}/J \neq 0/J$. Then $\alpha_{v_{1}}\cdots \alpha_{v_{j}}/J \neq 0/J$. Thus it follows immediately from Claim \ref{claim1}  that $\{u_{1}, \ldots, u_{i}\}, \{v_{1}, \ldots, v_{j}\} \in E'$,
and so $i=j=3$, or $i= j= 2$ with $\{u_{1}, u_{2}\}$ and $\{v_{1}, v_{2}\}$ being linked. This implies by Lemma \ref{lem1} that
${\bf a}_{u_{1}}\cdots{\bf a}_{u_{i}}/R^{\#} = {\bf a}_{v_{1}}\cdots{\bf a}_{v_{j}}/R^{\#}$, as required.

Conversely, assume that ${\bf a}_{u_{1}}\cdots{\bf a}_{u_{i}}/R^{\#} = {\bf a}_{v_{1}}\cdots{\bf a}_{v_{j}}/R^{\#}$. Then we can similarly get that  $\alpha_{u_{1}}\cdots \alpha_{u_{i}}/J=\alpha_{v_{1}}\cdots \alpha_{v_{j}}/J$. Thus we have
$$\alpha_{u_{1}}\cdots \alpha_{u_{i}}/J=\alpha_{v_{1}}\cdots \alpha_{v_{j}}/J \Longleftrightarrow {\bf a}_{u_{1}}\cdots{\bf a}_{u_{i}}/R^{\#}= {\bf a}_{v_{1}}\cdots{\bf a}_{v_{j}}/R^{\#}.$$

We have shown that $\psi$ is well-defined and a monomorphism. It also is easy to see that  $\psi$ is also an epimorphism.
Thus $\psi$ is an isomorphism and so $S_{\scriptscriptstyle \mathbb{H}}\in \mathbf{V}(S_c(abc))$.
\epf
As a consequence, we have immediately
\begin{cor}
If 3-hypergraph $\, \mathbb{H}$ is a 3-hyperforest or $n$-cycle except 3-cycle, then $\mathbf{V}(S_{\scriptscriptstyle \mathbb{H}})= \mathbf{V}(S_c(abc))$.
\end{cor}
For a 3-cycle 3-hypergraph, we also have
\begin{prop}\label{lem4}
Let $\mathbb{H}$ be a 3-cycle 3-hypergraph. Then $\mathbf{V}(S_c(abc)) \subset \mathbf{V}(S_{\scriptscriptstyle \mathbb{H}}) \subset \mathbf{V}(S_c(abcd))$.
\end{prop}
\pf
It is easy to see that $S_c(abc)\in \mathbf{V}(S_{\scriptscriptstyle \mathbb{H}})$ and so  $\mathbf{V} (S_c(abc)) \subseteq \mathbf{V}(S_{\scriptscriptstyle \mathbb{H}})$.

Also,  for any  substitution from $\{x_1, x_2, \ldots, x_6\}$ to $S_c(abc)$ we have
$$\varphi(x_1x_2x_3+x_3x_4x_5+x_5x_6x_1),  \varphi((x_1+x_4)(x_2+x_5)(x_3+x_6)) \in \{0, abc\}, \, and $$
$$\varphi(x_1x_2x_3+x_3x_4x_5+x_5x_6x_1) = abc \Longleftrightarrow   \varphi((x_1+x_4)(x_2+x_5)(x_3+x_6)) = abc.$$

This implies that $S_c(abc)$ satisfies the identity
\begin{equation} \label{f31}
x_1x_2x_3+x_3x_4x_5+x_5x_6x_1\approx (x_1+x_4)(x_2+x_5)(x_3+x_6).
\end{equation}
Assume that 3-cycle 3-hypergraph $\mathbb{H}$ is a  sequence $u_1, e_1, u_3, e_2, u_5, e_3$ , where
$$e_1 = \{u_1, u_2, u_3\}, e_2 = \{u_3, u_4, u_5\}, e_3 = \{u_5, u_6, u_1\}.$$
For the  substitution $\varphi$ from $\{x_1, x_2, \ldots, x_6\}$ to $S_{\scriptscriptstyle \mathbb{H}}$ defined by $\varphi(x_i) = \bar{{\bf a}}_{u_i}$,  we have
$$\varphi(x_1x_2x_3+x_3x_4x_5+x_5x_6x_1) = \bar{{\bf a}} \,\,but \,\,   \varphi((x_1+x_4)(x_2+x_5)(x_3+x_6)) = 0.$$
This implies that $S_{\scriptscriptstyle \mathbb{H}}$ does not satisfy the identity (\ref{f31}). This shows that
$$\mathbf{V}(S_c(abc))\subset \mathbf{V}(S_{\scriptscriptstyle \mathbb{H}}).$$

Let $(S_c(abcd))^2$ denote the direct product of 2 copies of $S_c(abcd)$. Consider the subsemiring $A$ of $(S_c(abcd))^2$ generated by $\{\alpha_1, \alpha_2, \ldots , \alpha_6\}$, where
$$\alpha_1=(a, bc), \alpha_2=(bc, d), \alpha_3=(d, a), \alpha_4=(ab, bc), \alpha_5=(c, d), \alpha_6=(bd, a).$$
Take
$$J=\{(x_1, x_2)\in A\mid (\exists~ t\in \{1, 2\})~ x_t=0 ~\text{or}~ (\exists~ t\in \{1, 2\})~ x_t=abcd, x_t\neq x_j ~(t\neq j)\}.$$
It is easy to verify that $J$ is an ideal of $A$ and $A/J$ is isomorphic to $S_{\scriptscriptstyle \mathbb{H}}$, under the obvious map defined on generators by $\bar{\alpha}_i\mapsto \bar{{\bf a}}_{u_i}$. This implies that $S_{\scriptscriptstyle \mathbb{H}}\in \mathbf{V}(S_c(abcd))$ and so $\mathbf{V}(S_{\scriptscriptstyle \mathbb{H}}) \subseteq \mathbf{V}(S_c(abcd))$.
Also, $S_{\scriptscriptstyle \mathbb{H}}$ satisfies $x_1x_2x_3x_4 \approx y_1y_2y_3y_4$ but $S_c(abcd)$ does not. Thus $\mathbf{V}(S_{\scriptscriptstyle \mathbb{H}})\subset \mathbf{V}(S_c(abcd))$.
Summarizing the above results we have
$$\mathbf{V}(S_c(abc)) \subset \mathbf{V}(S_{\scriptscriptstyle \mathbb{H}}) \subset \mathbf{V}(S_c(abcd)),$$
as required.
\epf

\section{ Varieties generated by 3-hypergraph semirings defined from some special 3-hypergraphs}

In this section we shall introduce and study 3-hypergraphs of a few other types , such as  beam-type hypergraph, fan-type hypergraph and nested-type hypergraph etc, and
varieties generated by their corresponding hypergraph semirings.

The following Lemma \ref{prop3.2} and Proposition \ref{cor3.3}  will tell us that in the case of girth being finite, it is only needed to study the varieties generated by the hypergraph semirings defined from 3-hypergraphs without leaf.

\begin{lem}\label{prop3.2}
	For any 3-hypergraph $\mathbb{H}$ except 3-hyperforest,
	there exists a 3-hypergraph $\mathbb{H}'$ without leaf such that
	$\mathbf{V}(S_{\scriptscriptstyle \mathbb{H}})=\mathbf{V}(S_{{\scriptscriptstyle \mathbb{H}}'})$.
\end{lem}
\pf
Let $\mathbb{H}=(V, E)$ with $V=\{u_1, \ldots, u_n, u_{n+1}, u_{n+2}, u_{n+3}\}$ and $E=\{e_1, e_2, \ldots, e_{m+1} \}$.
If $\mathbb{H}$ has no leaf, then the statement holds.
Otherwise, assume without loss of generality
that $e_{m+1}=\{u_{n+1}, u_{n+2}, u_{n+3}\}$ is a leaf in $\mathbb{H}$.
Notice that $\mathbb{H}$ has cycles, since it is not a hyperforest. Thus there exists a hyperedge different from $e_{m+1}$
in $\mathbb{H}$. Write it as $e_1=\{u_1, u_2, u_3\}$.

If $e_{m+1}\cap e_i=\emptyset$ for any $1\leq i \leq m$, then $\mathbf{V}(S_{\scriptscriptstyle \mathbb{H}})= \mathbf{V}(S_{{\scriptscriptstyle \mathbb{H}}_1})$,
where $\mathbb{H}_1$ denote the subhypergraph of $\mathbb{H}$ induced by the vertex set $V\setminus e_{m+1}$. In fact,
the hypergraph semiring $S_{{\scriptscriptstyle \mathbb{H}}_1}$ is generated by $\{\bar{{\bf a}}_{u_1}, \bar{{\bf a}}_{u_2}, \ldots, \bar{{\bf a}}_{u_n} \}$.
Clearly, $S_{{\scriptscriptstyle \mathbb{H}}_1}\in \mathbf{V}(S_{\scriptscriptstyle \mathbb{H}})$. We shall show that $S_{\scriptscriptstyle \mathbb{H}}\in \mathbf{V}(S_{{\scriptscriptstyle \mathbb{H}}_1})$.
Let $(S_{{\scriptscriptstyle \mathbb{H}}_1})^3$ denote the direct product of 3 copies of $S_{{\scriptscriptstyle \mathbb{H}}_1}$. Consider the subsemiring $A$ of $(S_{{\scriptscriptstyle \mathbb{H}}_1})^3$ generated by $$\{\alpha_1, \alpha_2, \ldots, \alpha_n , \alpha_{n+1}, \alpha_{n+2}, \alpha_{n+3}\},$$ where
$$\alpha_1=(\bar{{\bf a}}_{u_1}, \bar{{\bf a}}_{u_1}, \bar{{\bf a}}_{u_1}), \alpha_2=(\bar{{\bf a}}_{u_2}, \bar{{\bf a}}_{u_2}, \bar{{\bf a}}_{u_2}),
\alpha_3=(\bar{{\bf a}}_{u_3}, \bar{{\bf a}}_{u_3}, \bar{{\bf a}}_{u_3}),
\ldots, \alpha_n=(\bar{{\bf a}}_{u_n}, \bar{{\bf a}}_{u_n}, \bar{{\bf a}}_{u_n}), $$
$$\alpha_{n+1}=(\bar{{\bf a}}_{u_1}, \bar{{\bf a}}_{u_2}, \bar{{\bf a}}_{u_3}), \alpha_{n+2}=(\bar{{\bf a}}_{u_2}, \bar{{\bf a}}_{u_3}, \bar{{\bf a}}_{u_1}), \alpha_{n+3}=(\bar{{\bf a}}_{u_3}, \bar{{\bf a}}_{u_1}, \bar{{\bf a}}_{u_2}).$$
Take
$$J=\{(x_1, x_2, x_3)\in A\mid (\exists t\in \{1, 2, 3\})~ x_t=0 \}.$$
It is easy to verify that $J$ is an ideal of $A$ and $A/J$ is isomorphic to $S_{\scriptscriptstyle \mathbb{H}}$ under the obvious map defined on generators by $\bar{\alpha}_i\mapsto \bar{{\bf a}}_{u_i}$. Thus $S_{\scriptscriptstyle \mathbb{H}}\in \mathbf{V}(S_{{\scriptscriptstyle \mathbb{H}}_1})$ and so $\mathbf{V}(S_{\scriptscriptstyle \mathbb{H}})= \mathbf{V}(S_{{\scriptscriptstyle \mathbb{H}}_1})$

If $ e_{m+1}\cap e_i \neq \emptyset$ for some $1\leq i \leq m$, then $|e_{m+1}\cap e_i|=1$,  since $\mathbb{H}$ is linear.
Without loss of generality, we suppose that $|e_{m+1}\cap e_1|=1$ and
$u_1=u_{n+3} \in e_{m+1}\cap e_1$. Let $\mathbb{H}_1=(V_1, E_1)$ denote the partial subhypergraph of $\mathbb{H}$ where $E_1=E\setminus \{e_{m+1}\}$ and $V_1=\cup\{e_i\mid e_i\in E_1\}$.
Thus we can construct a hypergraph semiring $S_{{\scriptscriptstyle \mathbb{H}}_1}$, which is generated by $\{\bar{{\bf a}}_{u_1}, \bar{{\bf a}}_{u_2}, \ldots, \bar{{\bf a}}_{u_n} \}$.
Clearly, $S_{{\scriptscriptstyle \mathbb{H}}_1}\in \mathbf{V}(S_{\scriptscriptstyle \mathbb{H}})$. We shall prove that $S_{\scriptscriptstyle \mathbb{H}}\in \mathbf{V}(S_{{\scriptscriptstyle \mathbb{H}}_1})$.
Let $(S_{{\scriptscriptstyle \mathbb{H}}_1})^3$ denote the direct product of 3 copies of $S_{{\scriptscriptstyle \mathbb{H}}_1}$. Consider the subsemiring $A$ of $(S_{{\scriptscriptstyle \mathbb{H}}_1})^3$ generated by $$\{\alpha_1, \alpha_2, \ldots, \alpha_n , \alpha_{n+1}, \alpha_{n+2}\},$$ where
$$\alpha_1=(\bar{{\bf a}}_{u_1}, \bar{{\bf a}}_{u_1}, \bar{{\bf a}}_{u_1}), \alpha_2=(\bar{{\bf a}}_{u_2}, \bar{{\bf a}}_{u_2}, \bar{{\bf a}}_{u_2}),
\alpha_3=(\bar{{\bf a}}_{u_3}, \bar{{\bf a}}_{u_3}, \bar{{\bf a}}_{u_3}),
\ldots, \alpha_n=(\bar{{\bf a}}_{u_n}, \bar{{\bf a}}_{u_n}, \bar{{\bf a}}_{u_n}), $$
$$\alpha_{n+1}=(\bar{{\bf a}}_{u_2}, \bar{{\bf a}}_{u_3}, \bar{{\bf a}}_{u_2}), \alpha_{n+2}=(\bar{{\bf a}}_{u_3}, \bar{{\bf a}}_{u_2}, \bar{{\bf a}}_{u_3}).$$
Take
$$J=\{(x_1, x_2, x_3)\in A\mid (\exists t\in \{1, 2, 3\})~ x_t=0 \}.$$
It is easy to verify that $J$ is an ideal of $A$ and $A/J$ is isomorphic to $S_{\scriptscriptstyle \mathbb{H}}$ under the obvious map defined on generators by $\bar{\alpha}_i\mapsto \bar{{\bf a}}_{u_i}$. Thus $S_{\scriptscriptstyle \mathbb{H}}\in \mathbf{V}(S_{{\scriptscriptstyle \mathbb{H}}_1})$ and so $\mathbf{V}(S_{\scriptscriptstyle \mathbb{H}})= \mathbf{V}(S_{{\scriptscriptstyle \mathbb{H}}_1})$.

If $\mathbb{H}_1$ has no leaf, then we may set $\mathbb{H}' = \mathbb{H}_1$, as required.
Otherwise, suppose that $\mathbb{H}_1=(V_1, E_1)$ has a leaf $e_{m}$. Write  $\mathbb{H}_2=(V_2, E_2)$ to denote the partial subhypergraph of $\mathbb{H}_1$
induced by the hyperedge set $E_2=E_1\setminus \{e_{m}\}$ and the vertex set $V_2=\cup\{e_i\mid e_i\in E_2\}$.
Similarly, we can get that $\mathbf{V}(S_{{\scriptscriptstyle \mathbb{H}}_1}) = \mathbf{V}(S_{{\scriptscriptstyle \mathbb{H}}_2})$ and so  $\mathbf{V}(S_{\scriptscriptstyle \mathbb{H}}) = \mathbf{V}(S_{{\scriptscriptstyle \mathbb{H}}_2})$. If $\mathbb{H}_2$ has no leaf, then we may set $\mathbb{H}' = \mathbb{H}_2$, as required. Otherwise, we are going on the above procedure. At most steps, we can see that
there exists a hypergraph $\mathbb{H}_n = (V_n, E_n)$ without leaf
(i.e., the partial subhypergraph of $\mathbb{H}$ induced by $E_n=\{e\in E \mid \text{$e$ is not a leaf} \}$ and $V_n=\cup\{e \mid e\in E'\}$)
such that $\mathbf{V}(S_{\scriptscriptstyle \mathbb{H}})=\mathbf{V}(S_{{\scriptscriptstyle \mathbb{H}}_{n}})$.
Thus we may take $\mathbb{H}' = \mathbb{H}_n$. This completes our proof.
\epf

\begin{prop}\label{cor3.3}
	Each variety generated by 3-hypergraph semirings except 3-hyperforest semirings is equal to a variety generated by 3-hypergraph semirings, which hypergraphs all have no leaf.
\end{prop}
\pf
It follows directly from Lemma \ref{prop3.2} that for any given the family $\{\mathbb{H}_i \mid i \in I\}$ of 3-hypergraphs,
there exists a family $\{\mathbb{H}^\prime_i \mid i \in I\}$ of 3-hypergraphs without leaf  such that  $ \mathbf{V}(S_{{\scriptscriptstyle \mathbb{H}}_i}) = \mathbf{V}(S_{{\scriptscriptstyle \mathbb{H}}^\prime_i})$ for all $i \in I$.
This implies that
\begin{equation*}
\mathbf{V}({\displaystyle\bigcup}_{i \in I} S_{{\scriptscriptstyle \mathbb{H}}_i})
= {\displaystyle \bigvee}_{i \in I} \mathbf{V}(S_{{\scriptscriptstyle \mathbb{H}}_i})
={\displaystyle \bigvee}_{i \in I} \mathbf{V}(S_{{\scriptscriptstyle \mathbb{H}}^\prime_i})
=\mathbf{V}({\displaystyle\bigcup}_{i \in I} S_{{\scriptscriptstyle \mathbb{H}}^\prime_i}).  \qedhere
\end{equation*}
\epf

\begin{figure}[H]
	\centering
	\captionsetup[subfloat]{labelsep=none, format=plain, labelformat=empty}		
	\subfloat[$\mathbb{B}_1$]{	
		\scalebox{0.7}{	
			\begin{tikzpicture}[style=thick]
			\draw (0,0) -- (1,0) -- (2,0) -- (1.5,0.8660254037844386) -- (1,1.732050807568877) -- (0.5,0.8660254037844386) -- cycle;
			\fill (0,0)  circle (2pt) ;
			\fill (1,0)  circle (2pt) ;
			\fill (2,0)  circle (2pt) ;
			\fill (1.5,0.8660254037844386)   circle  (2pt) ;
			\fill (1,1.732050807568877)  circle  (2pt) ;
			\fill (0.5,0.8660254037844386)  circle  (2pt) ;
			\end{tikzpicture}
	}}
	\subfloat[$\mathbb{B}_2$]{
		\scalebox{0.7}{	
			\begin{tikzpicture}[style=thick]
			\draw (0,0) -- (1,0) -- (2,0) -- (2.5,0.8660254037844386) -- (3,1.732050807568877) -- (2,1.732050807568877) -- (1,1.732050807568877) --(0.5,0.8660254037844386) -- cycle;
			\draw (1,1.732050807568877) -- (1.5,0.8660254037844386) -- (2,0);
			\fill (0,0)  circle  (2pt) ;
			\fill (1,0)  circle  (2pt) ;
			\fill (2,0)  circle  (2pt) ;
			\fill (2.5,0.8660254037844386)  circle  (2pt) ;
			\fill (3,1.732050807568877)  circle  (2pt) ;
			\fill (2,1.732050807568877)  circle  (2pt) ;
			\fill (1,1.732050807568877)  circle  (2pt) ;
			\fill (0.5,0.8660254037844386)  circle  (2pt) ;
			\fill (1.5,0.8660254037844386) circle  (2pt) ;
			\end{tikzpicture}
	}}	
	\subfloat[$\mathbb{B}_3$]{
		\scalebox{0.7}{	
			\begin{tikzpicture}[style=thick]
			\draw (0,0) -- (1,0) -- (2,0) -- (3,0) -- (4,0) -- (3.5,0.8660254037844386) -- (3,1.732050807568877) -- (2,1.732050807568877) -- (1,1.732050807568877) --(0.5,0.8660254037844386) -- cycle;
			\draw (1,1.732050807568877) -- (1.5,0.8660254037844386) -- (2,0);
			\draw (3,1.732050807568877) -- (2.5,0.8660254037844386) -- (2,0);
			\fill (0,0)  circle  (2pt) ;
			\fill (1,0)  circle  (2pt) ;
			\fill (2,0) circle  (2pt) ;
			\fill (3,0) circle  (2pt) ;
			\fill (4,0) circle  (2pt) ; 	
			\fill (3.5,0.8660254037844386)  circle  (2pt) ;
			\fill (3, 1.732050807568877)  circle  (2pt) ;
			\fill (2,1.732050807568877)  circle  (2pt) ;
			\fill (1,1.732050807568877) circle  (2pt) ;
			\fill (0.5,0.8660254037844386)  circle  (2pt) ;
			\fill (1.5,0.8660254037844386)  circle  (2pt) ;
			\fill (2.5,0.8660254037844386) circle  (2pt) ;
			\end{tikzpicture}
	}}	
	\subfloat[$\mathbb{B}_4$]{
		\scalebox{0.7}{	
			\begin{tikzpicture}[style=thick]
			\draw (0,0) -- (1,0) -- (2,0) -- (3,0) -- (4,0) -- (4.5, 0.8660254037844386)-- (5, 1.732050807568877) -- (4 ,1.732050807568877) -- (3,1.732050807568877) -- (2,1.732050807568877) -- (1,1.732050807568877) --(0.5,0.8660254037844386) -- cycle;
			\draw (1,1.732050807568877) -- (1.5,0.8660254037844386) -- (2,0);
			\draw (3,1.732050807568877) -- (2.5,0.8660254037844386) -- (2,0);
			\draw (3,1.732050807568877) -- (3.5,0.8660254037844386) -- (4,0);
			\fill (0,0)  circle  (2pt) ;
			\fill (1,0)  circle  (2pt) ;
			\fill (2,0)  circle  (2pt) ;
			\fill (3,0)  circle  (2pt) ;
			\fill (4,0)  circle  (2pt) ; 	
			\fill (4.5,0.8660254037844386)  circle  (2pt) ;
			\fill (5, 1.732050807568877)  circle  (2pt) ;
			\fill (4,1.732050807568877)  circle  (2pt) ;
			\fill (3,1.732050807568877)  circle  (2pt) ;
			\fill (2,1.732050807568877)  circle  (2pt) ;
			\fill (1,1.732050807568877)  circle  (2pt) ;
			\fill (0.5,0.8660254037844386)  circle  (2pt) ;
			\fill (1.5,0.8660254037844386)  circle  (2pt) ;
			\fill (2.5,0.8660254037844386)  circle  (2pt) ;
			\fill (3.5,0.8660254037844386)  circle  (2pt) ;
			\end{tikzpicture}
	}}
	\hspace{5mm} \ldots \hspace{5mm}
	\caption{Beam-type 3-hypergraph}\label{tu1}
\end{figure}
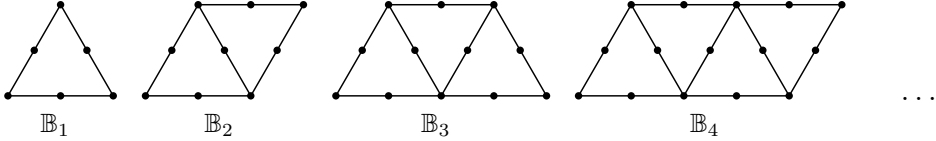

In the following we investigate varieties generated by a 3-hypergraph semiring $S_{\scriptscriptstyle \mathbb{B}_i}$, where each 3-hypergraph $\mathbb{B}_i$ is given in figure \ref{tu1}, and  all hyperedges in each $\mathbb{B}_i$ are only denoted by straight lines and all vertices are denoted by black dots.
\vspace{-3mm}
\vspace{-3mm}
\noindent We refer to $\mathbb{B}_i$ as the $i$-th beam-type 3-hypergraph. In particular, We also refer to the 3-cycle 3-hypergraph shaped like $\mathbb{B}_1$ as a triangle.

\begin{prop}\label{prop3.4}
	For all $i\geq 1$,
	$\mathbf{V}(S_{{\scriptscriptstyle \mathbb{B}}_{i}})= \mathbf{V}(S_{{\scriptscriptstyle \mathbb{B}}_{1}})$.
\end{prop}
\pf We shall only prove that $\mathbf{V}(S_{{\scriptscriptstyle \mathbb{B}}_{i}})= \mathbf{V}(S_{{\scriptscriptstyle \mathbb{B}}_{i+1}})$ for all $i\geq 1$.
It is easy to see that $\mathbf{V}(S_{{\scriptscriptstyle \mathbb{B}}_{i}})\subseteq \mathbf{V}(S_{{\scriptscriptstyle \mathbb{B}}_{i+1}})$, since $S_{{\scriptscriptstyle \mathbb{B}}_{i}}$
is a subsemiring of $S_{{\scriptscriptstyle \mathbb{B}}_{i+1}}$ up to isomorphism.

To show that $\mathbf{V}(S_{{\scriptscriptstyle \mathbb{B}}_{i+1}})\subseteq \mathbf{V}(S_{{\scriptscriptstyle \mathbb{B}}_{i}})$, We need only to prove that
$S_{{\scriptscriptstyle \mathbb{B}}_{i+1}}\in \mathbf{V}(S_{{\scriptscriptstyle \mathbb{B}}_{i}})$ for all $i \geq 1$.
We shall firstly prove that $S_{{\scriptscriptstyle \mathbb{B}}_{2}}\in \mathbf{V}(S_{{\scriptscriptstyle \mathbb{B}}_{1}})$.
If the vertices of $\mathbb{B}_1$ and $\mathbb{B}_2$ are labelled as follows
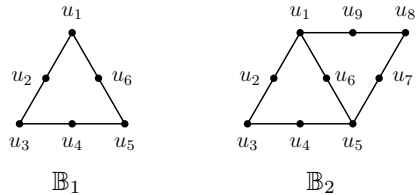
\vspace{-3mm}
\begin{figure}[H]
	\centering  	
	\captionsetup[subfloat]{labelsep=none, format=plain, labelformat=empty}	
	\subfloat[$\mathbb{B}_1$]{
		\scalebox{0.7}{		
			\begin{tikzpicture}[style=thick]
			\draw (0,0) -- (1,0) -- (2,0) -- (1.5,0.8660254037844386) -- (1,1.732050807568877) -- (0.5,0.8660254037844386) -- cycle;
			\fill (0,0) node[below=1mm]{$u_3$} circle (2pt) ;
			\fill (1,0) node[below=1mm]{$u_4$} circle (2pt) ;
			\fill (2,0) node[below=1mm]{$u_5$} circle (2pt) ;
			\fill (1.5,0.8660254037844386) node[right=1mm]{$u_6$}  circle  (2pt) ;
			\fill (1,1.732050807568877) node[above=1mm]{$u_1$} circle  (2pt) ;
			\fill (0.5,0.8660254037844386) node[left=1mm]{$u_2$} circle  (2pt) ;
			\end{tikzpicture}
	}}
	\qquad
	\subfloat[$\mathbb{B}_2$]{
		\scalebox{0.7}{	
			\begin{tikzpicture}[style=thick]
			\draw (0,0) -- (1,0) -- (2,0) -- (2.5,0.8660254037844386) -- (3,1.732050807568877) -- (2,1.732050807568877) -- (1,1.732050807568877) --(0.5,0.8660254037844386) -- cycle;
			\draw (1,1.732050807568877) -- (1.5,0.8660254037844386) -- (2,0);
			\fill (0,0) node[below=1mm]{$u_3$} circle (2pt) ;
			\fill (1,0) node[below=1mm]{$u_4$} circle (2pt) ;
			\fill (2,0) node[below=1mm]{$u_5$} circle (2pt) ;
			\fill (2.5,0.8660254037844386) node[ right=1mm]{$u_7$} circle  (2pt) ;
			\fill (3,1.732050807568877) node[above=1mm]{$u_8$} circle  (2pt) ;
			\fill (2,1.732050807568877) node[above=1mm]{$u_9$} circle  (2pt) ;
			\fill (1,1.732050807568877) node[above=1mm]{$u_1$} circle  (2pt) ;
			\fill (0.5,0.8660254037844386) node[ left=1mm]{$u_2$} circle  (2pt) ;
			\fill (1.5,0.8660254037844386) node[right=0.1mm]{$u_6$} circle  (2pt) ;
			\end{tikzpicture}
	}}	
\end{figure}
\vspace{-3mm}
\noindent then it is well-known that hypergraph semiring $S_{{\scriptscriptstyle \mathbb{B}}_1}$ is generated by $\{\bar{{\bf a}}_{u_1}, \bar{{\bf a}}_{u_2}, \ldots, \bar{{\bf a}}_{u_6} \}$.
Write $(S_{{\scriptscriptstyle \mathbb{B}}_1})^3$ to denote the direct product of 3 copies of $S_{{\scriptscriptstyle \mathbb{B}}_1}$ and $A_1$ to denote the subsemiring of $(S_{{\scriptscriptstyle \mathbb{B}}_1})^3$ generated by $\{ \alpha_1, \ldots, \alpha_6, \beta_{11}, \beta_{12}, \beta_{13} \},$ where
$$\alpha_1 =(\bar{{\bf a}}_{u_1},  \bar{{\bf a}}_{u_1}, \bar{{\bf a}}_{u_1}), \, \alpha_2=(\bar{{\bf a}}_{u_2}, \bar{{\bf a}}_{u_2}, \bar{{\bf a}}_{u_3}),\,
\alpha_3=(\bar{{\bf a}}_{u_3},  \bar{{\bf a}}_{u_3},  \bar{{\bf a}}_{u_2}),$$
$$ \alpha_4 =(\bar{{\bf a}}_{u_4},  \bar{{\bf a}}_{u_4}, \bar{{\bf a}}_{u_1}), \, \alpha_5=(\bar{{\bf a}}_{u_5}, \bar{{\bf a}}_{u_5}, \bar{{\bf a}}_{u_3}), \, \alpha_6=(\bar{{\bf a}}_{u_6}, \bar{{\bf a}}_{u_6}, \bar{{\bf a}}_{u_2}),$$
$$\beta_{11}=(\bar{{\bf a}}_{u_4}, \bar{{\bf a}}_{u_1}, \bar{{\bf a}}_{u_4}), \, \beta_{12}=(\bar{{\bf a}}_{u_3}, \bar{{\bf a}}_{u_6}, \bar{{\bf a}}_{u_5}), \,\beta_{13}=(\bar{{\bf a}}_{u_2}, \bar{{\bf a}}_{u_5}, \bar{{\bf a}}_{u_6}).$$
Take
$$J_1=\{(x_1, x_2, x_3)\in A_1\mid (\exists \,t\in \{1, 2, 3\})~ x_t=0 \}.$$
It is easy to see that $J_1$ is an ideal of $A_1$ and factor semiring $A_1/J_1$ is generated by $\{\bar{\alpha}_1, \ldots, \bar{\alpha}_6, \bar{\beta}_{11}, \bar{\beta}_{12}, \bar{\beta}_{13}\}$.
Also, it is easy to verify that $A_1/J_1$ is isomorphic to hypergraph semiring $S_{\overline{\scriptscriptstyle \mathbb{B}}_{2}}$, where hypergraph $\overline{\mathbb{B}}_{2}=(\overline{V}_{2}, \overline{E}_{2})$ with $\overline{V}_{2}=\{\bar{\alpha}_1, \ldots, \bar{\alpha}_6, \bar{\beta}_{11}, \bar{\beta}_{12}, \bar{\beta}_{13}\}$
and
$$\overline{E}_2=\{\{\bar{\alpha}_1,\bar{\alpha}_2,\bar{\alpha}_3\}, \{\bar{\alpha}_3,\bar{\alpha}_4,\bar{\alpha}_5\}, \{\bar{\alpha}_5,\bar{\alpha}_6,\bar{\alpha}_1\}, \{\bar{\alpha}_5, \bar{\beta}_{11}, \bar{\beta}_{12}\},  \{\bar{\beta}_{12}, \bar{\beta}_{13}, \bar{\alpha}_1\} \}.$$
On the other hand, the map $\varphi_{2}: \overline{V}_{2}\rightarrow {V}_{2}$ defined by
$$\varphi_{2}(\bar{\alpha}_i ) = u_i \,\, (1 \leq i \leq 6), \,\, \varphi_{2}(\bar{\beta}_{11} ) = u_{7}, \,\,
\varphi_{2}(\bar{\beta}_{12} ) = u_{8}\,\, and \,\,\varphi_{2}(\bar{\beta}_{13} ) = u_{9}$$
is a hypergraph isomorphism from  $\overline{\mathbb{B}}_{2}$ to $\mathbb{B}_{2}$ and so $S_{\overline{\scriptscriptstyle \mathbb{B}}_{2}} \cong S_{{ \scriptscriptstyle \mathbb{B}}_{2}}$.
This implies that  $A_{1}/J_{1}\cong S_{{\scriptscriptstyle \mathbb{B}}_{2}}$ and so
$ S_{{\scriptscriptstyle \mathbb{B}}_{2}} \in  \mathbf{V}(S_{{\scriptscriptstyle \mathbb{B}}_{1}})$, as required.

To prove that $S_{{\scriptscriptstyle \mathbb{B}}_{3}}\in \mathbf{V}(S_{{\scriptscriptstyle \mathbb{B}}_{2}})$,
we label the vertices of $\mathbb{B}_2$ and $\mathbb{B}_3$ as follows
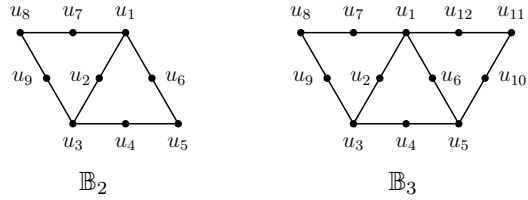
\vspace{-3mm}
\begin{figure}[H]
	\centering  	
	\captionsetup[subfloat]{labelsep=none, format=plain, labelformat=empty}	
	\subfloat[$\mathbb{B}_2$]{
		\scalebox{0.7}{		
			\begin{tikzpicture}[style=thick]
			\draw (0,0) -- (1,0) -- (2,0) -- (1.5,0.8660254037844386) -- (1,1.732050807568877) -- (0.5,0.8660254037844386) -- cycle;
			\draw (0,0) --(-0.5,0.8660254037844386)--(-1,1.732050807568877)--(0,1.732050807568877)--(1,1.732050807568877);
			\fill (0,0) node[below=1mm]{$u_3$} circle (2pt) ;
			\fill (1,0) node[below=1mm]{$u_4$} circle (2pt) ;
			\fill (2,0) node[below=1mm]{$u_5$} circle (2pt) ;
			\fill (1.5,0.8660254037844386) node[right=1mm]{$u_6$}  circle  (2pt) ;
			\fill (1,1.732050807568877) node[above=1mm]{$u_1$} circle  (2pt) ;
			\fill (0.5,0.8660254037844386) node[left=0.1mm]{$u_2$} circle  (2pt) ;
			\fill (-0.5,0.8660254037844386) node[left=1mm]{$u_9$} circle  (2pt) ;
			\fill (-1,1.732050807568877) node[above=1mm]{$u_8$} circle  (2pt) ;
			\fill (0,1.732050807568877) node[above=1mm]{$u_7$} circle  (2pt) ;
			\end{tikzpicture}
	}}
	\qquad
	\subfloat[$\mathbb{B}_3$]{
		\scalebox{0.7}{	
			\begin{tikzpicture}[style=thick]
			\draw (0,0) -- (1,0) -- (2,0) -- (2.5,0.8660254037844386) -- (3,1.732050807568877) -- (2,1.732050807568877) -- (1,1.732050807568877) --(0.5,0.8660254037844386) -- cycle;
			\draw (1,1.732050807568877) -- (1.5,0.8660254037844386) -- (2,0);
			\draw (0,0) --(-0.5,0.8660254037844386)--(-1,1.732050807568877)--(0,1.732050807568877)--(1,1.732050807568877);
			\fill (0,0) node[below=1mm]{$u_3$} circle (2pt) ;
			\fill (1,0) node[below=1mm]{$u_4$} circle (2pt) ;
			\fill (2,0) node[below=1mm]{$u_5$} circle (2pt) ;
			\fill (2.5,0.8660254037844386) node[ right=1mm]{$u_{10}$} circle  (2pt) ;
			\fill (3,1.732050807568877) node[above=1mm]{$u_{11}$} circle  (2pt) ;
			\fill (2,1.732050807568877) node[above=1mm]{$u_{12}$} circle  (2pt) ;
			\fill (1,1.732050807568877) node[above=1mm]{$u_1$} circle  (2pt) ;
			\fill (0.5,0.8660254037844386) node[ left=0.1mm]{$u_2$} circle  (2pt) ;
			\fill (1.5,0.8660254037844386) node[right=0.1mm]{$u_6$} circle  (2pt) ;
			\fill (-0.5,0.8660254037844386) node[left=1mm]{$u_9$} circle  (2pt) ;
			\fill (-1,1.732050807568877) node[above=1mm]{$u_8$} circle  (2pt) ;
			\fill (0,1.732050807568877) node[above=1mm]{$u_7$} circle  (2pt) ;
			\end{tikzpicture}
	}}
\end{figure}
\vspace{-3mm}
\noindent Thus it is well-known  that hypergraph semiring $S_{{\scriptscriptstyle \mathbb{B}}_2}$ is generated by $\{\bar{{\bf a}}_{u_1}, \bar{{\bf a}}_{u_2}, \ldots, \bar{{\bf a}}_{u_9} \}$.
Write $(S_{{\scriptscriptstyle \mathbb{B}}_2})^3$ to denote the direct product of 3 copies of $S_{{\scriptscriptstyle \mathbb{B}}_2}$ and $A_2$ to denote
the subsemiring of $(S_{{\scriptscriptstyle \mathbb{B}}_2})^3$ generated by $\{\alpha_1, \ldots, \alpha_6, \beta_{11}, \beta_{12}, \beta_{13}, \beta_{21}, \beta_{22}, \beta_{23}\}$, where
$$\beta_{21}=(\bar{{\bf a}}_{u_7}, \bar{{\bf a}}_{u_7}, \bar{{\bf a}}_{u_2}), \,
\beta_{22}=(\bar{{\bf a}}_{u_8}, \bar{{\bf a}}_{u_8}, \bar{{\bf a}}_{u_3}),  \,
\beta_{23}=(\bar{{\bf a}}_{u_9}, \bar{{\bf a}}_{u_9}, \bar{{\bf a}}_{u_1}).$$
Take
$$J_2=\{(x_1, x_2, x_3)\in A_2\mid (\exists t\in \{1, 2, 3\})~ x_t=0 \}.$$
It is easy to verify that $J_2$ is an ideal of $A_2$ and factor semiring $A_2/J_2$ is generated by $\{\bar{\alpha}_1, \ldots, \bar{\alpha}_6, \bar{\beta}_{11}, \bar{\beta}_{12}, \bar{\beta}_{13}, \bar{\beta}_{21}, \bar{\beta}_{22}, \bar{\beta}_{23}\}$.
Also, it is easy to verify that $A_2/J_2$ is isomorphic to hypergraph semiring $S_{\overline{\scriptscriptstyle\mathbb{B}}_{3}}$, where hypergraph $\overline{\mathbb{B}}_{3}=(\overline{V}_{3}, \overline{E}_{3})$ with
$$\overline{V}_{3} = \overline{V}_{2} \cup \{\bar{\beta}_{21}, \bar{\beta}_{22}, \bar{\beta}_{23}\}  = \{\bar{\alpha}_1, \ldots, \bar{\alpha}_6, \bar{\beta}_{11}, \bar{\beta}_{12}, \bar{\beta}_{13}, \bar{\beta}_{21}, \bar{\beta}_{22}, \bar{\beta}_{23}\}$$
and $\overline{E}_{3}= \overline{E}_{2}\cup \{\{\bar{\alpha}_1, \bar{\beta}_{21}, \bar{\beta}_{22}\}, \{\bar{\beta}_{22}, \bar{\beta}_{23}, \bar{\alpha}_3 \}\}$.
On the other hand, define the map $\varphi_{3}: \overline{V}_{3}\rightarrow {V}_{3}$ by
$$\varphi_{3}(\bar{\alpha}_i ) = \varphi_{2}(\bar{\alpha}_i ) \,\, (1 \leq i \leq 6), \,\, \varphi_{3}(\bar{\beta}_{11} ) = u_{10}, \,\,
\varphi_{3}(\bar{\beta}_{12} ) = u_{11},\,\,\varphi_{3}(\bar{\beta}_{13} ) = u_{12},$$
$\varphi_{3}(\bar{\beta}_{21} ) = u_{7}, \,\,\varphi_{3}(\bar{\beta}_{22} ) = u_{8}$ and  $\varphi_{3}(\bar{\beta}_{23} ) = u_{9}$.
Then $\varphi_{3}$ is a hypergraph isomorphism  from $\overline{\mathbb{B}}_{3}$ to $\mathbb{B}_{3}$ and so $S_{\overline{\scriptscriptstyle\mathbb{B}}_{3}} \cong S_{{\scriptscriptstyle\mathbb{B}}_{3}}$.
This implies that  $A_{2}/J_{2}\cong S_{{\scriptscriptstyle\mathbb{B}}_{3}}$ and so $ S_{{\scriptscriptstyle\mathbb{B}}_{3}} \in  \mathbf{V}(S_{{\scriptscriptstyle \mathbb{B}}_{2}})$, as required.

To prove that $S_{{\scriptscriptstyle \mathbb{B}}_{i+1}} \in \mathbf{V}(S_{{\scriptscriptstyle \mathbb{B}}_{i}})$ ($i\geq 3$),
we label the vertices of $\mathbb{B}_i$ and $\mathbb{B}_{i+1}$ counterclockwise as follows
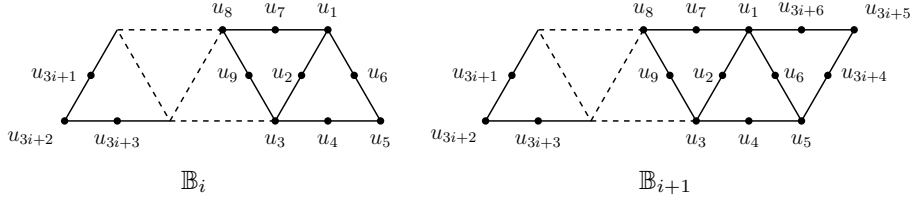
\vspace{-3mm}
\begin{figure}[H]
	\centering
	\captionsetup[subfloat]{labelsep=none, format=plain, labelformat=empty}	
	\subfloat[$\mathbb{B}_i$]{
		\scalebox{0.7}{	
			\begin{tikzpicture}[style=thick]
			\draw[dashed] (4,0) -- (3,0) -- (2,0) --(1.5,0.8660254037844386)--(1,1.732050807568877)-- (2,1.732050807568877) -- (3,1.732050807568877) ;
			\draw[dashed] (3,1.732050807568877) -- (2.5, 0.8660254037844386) -- (2,0) ;	
			\draw (4,0)--(3.5, 0.8660254037844386)--(3,1.732050807568877)--(4,1.732050807568877)--(5,1.732050807568877)--(5.5, 0.8660254037844386)--(6,0)--(5,0) -- cycle;
			\draw (5,1.732050807568877) -- (4.5,0.8660254037844386) --(4,0);
			\draw (0,0)--(1,0)--(2,0);
			\draw (1,1.732050807568877) -- (0.5,0.8660254037844386) --(0,0);
			\fill (0,0) node[below left=0.5mm]{$u_{3i+2}$} circle  (2pt) ;
			\fill (1,0) node[below=1mm]{$u_{3i+3}$} circle  (2pt) ;
			\fill (4,0) node[below=1mm]{$u_{3}$} circle  (2pt) ;
			\fill (5,0) node[below=1mm]{$u_{4}$} circle  (2pt) ;
			\fill (6,0)  node[below=1mm]{$u_{5}$} circle  (2pt) ; 	
			\fill (4.5,0.8660254037844386) node[left=0.1mm]{$u_{2}$} circle  (2pt) ;
			\fill (5.5,0.8660254037844386) node[right=1mm]{$u_{6}$} circle  (2pt) ;
			\fill (5,1.732050807568877) node[above=1mm]{$u_{1}$} circle  (2pt) ;
			\fill (0.5,0.8660254037844386) node[left=1mm]{$u_{3i+1}$} circle  (2pt) ;
			\fill (4,1.732050807568877) node[above=1mm]{$u_{7}$} circle  (2pt) ;
			\fill (3,1.732050807568877) node[above=1mm]{$u_{8}$} circle  (2pt) ;
			\fill (3.5,0.8660254037844386) node[left=0.5mm]{$u_{9}$} circle  (2pt) ;
			\end{tikzpicture}
	}}
	\subfloat[$\mathbb{B}_{i+1}$]{
		\scalebox{0.7}{	
			\begin{tikzpicture}[style=thick]
			\draw[dashed] (4,0) -- (3,0) -- (2,0) --(1.5,0.8660254037844386)--(1,1.732050807568877)-- (2,1.732050807568877) -- (3,1.732050807568877) ;
			\draw[dashed] (3,1.732050807568877) -- (2.5, 0.8660254037844386) -- (2,0) ;	
			\draw (4,0)--(3.5, 0.8660254037844386)--(3,1.732050807568877)--(4,1.732050807568877)--(5,1.732050807568877)--(5.5, 0.8660254037844386)--(6,0)--(5,0) -- cycle;
			\draw (5,1.732050807568877) -- (4.5,0.8660254037844386) --(4,0);
			\draw (6,0)--(6.5, 0.8660254037844386)--(7, 1.732050807568877)--(6, 1.732050807568877)--(5, 1.732050807568877);
			\draw (0,0)--(1,0)--(2,0);
			\draw (1,1.732050807568877) -- (0.5,0.8660254037844386) --(0,0);
			\fill (0,0) node[below left=0.1mm]{$u_{3i+2}$} circle  (2pt) ;
			\fill (1,0) node[below=1mm]{$u_{3i+3}$} circle  (2pt) ;
			\fill (4,0) node[below=1mm]{$u_{3}$} circle  (2pt) ;
			\fill (5,0) node[below=1mm]{$u_{4}$} circle  (2pt) ;
			\fill (6,0)  node[below=1mm]{$u_{5}$} circle  (2pt) ; 	
			\fill (4.5,0.8660254037844386) node[left=0.1mm]{$u_{2}$} circle  (2pt) ;
			\fill (5.5,0.8660254037844386) node[right=0.1mm]{$u_{6}$} circle  (2pt) ;
			\fill (5,1.732050807568877) node[above=1mm]{$u_{1}$} circle  (2pt) ;
			\fill (0.5,0.8660254037844386) node[left=1mm]{$u_{3i+1}$} circle  (2pt) ;
			\fill (4,1.732050807568877) node[above=1mm]{$u_{7}$} circle  (2pt) ;
			\fill (3,1.732050807568877) node[above=1mm]{$u_{8}$} circle  (2pt) ;
			\fill (3.5,0.8660254037844386) node[left=0.5mm]{$u_{9}$} circle  (2pt) ;
			\fill (6.5,0.8660254037844386) node[right=1mm]{$u_{3i+4}$} circle  (2pt) ;
			\fill (7,1.732050807568877) node[above right=0.5mm]{$u_{3i+5}$} circle  (2pt) ;
			\fill (6,1.732050807568877) node[above=1mm]{$u_{3i+6}$} circle  (2pt) ;
			\end{tikzpicture}
	}}
\end{figure}
\vspace{-3mm}
\noindent It is well-known  that hypergraph semiring $S_{{\scriptscriptstyle \mathbb{B}}_i}$ is generated by $\{\bar{{\bf a}}_{u_1}, \bar{{\bf a}}_{u_2}, \ldots, \bar{{\bf a}}_{u_{3i+3}} \}$.
Write $(S_{{\scriptscriptstyle \mathbb{B}}_{i}})^3$ to denote the direct product of 3 copies of $S_{{\scriptscriptstyle \mathbb{B}}_{i}}$  and $A_i$ to denote the subsemiring of $(S_{{\scriptscriptstyle \mathbb{B}}_i})^3$ generated by $\{\alpha_1, \ldots, \alpha_6, \beta_{11}, \beta_{12}, \beta_{13}, \ldots,  \beta_{(i-1)1}, \beta_{(i-1)2}, \beta_{(i-1)3}\}\cup \{\beta_{i1}, \beta_{i2}, \beta_{i3}\}$, where
$$\beta_{i1}=(\bar{{\bf a}}_{u_{3i+1}}, \bar{{\bf a}}_{u_{3i+1}}, x), \,
\beta_{i2}=(\bar{{\bf a}}_{u_{3i+2}}, \bar{{\bf a}}_{u_{3i+2}}, y),  \,
\beta_{i3}=(\bar{{\bf a}}_{u_{3i+3}}, \bar{{\bf a}}_{u_{3i+3}}, z),$$
$x,\, y$ and $z$ are equal to the $(1,k),\,(2,k),\,(3,k)$-entries of matrix
\[
\begin{pmatrix}
 \bar{{\bf a}}_{u_3} & \bar{{\bf a}}_{u_2} & \bar{{\bf a}}_{u_2} & \bar{{\bf a}}_{u_1} & \bar{{\bf a}}_{u_1} & \bar{{\bf a}}_{u_3}\\
\bar{{\bf a}}_{u_2} & \bar{{\bf a}}_{u_3} & \bar{{\bf a}}_{u_1} & \bar{{\bf a}}_{u_2} & \bar{{\bf a}}_{u_3} & \bar{{\bf a}}_{u_1}\\
\bar{{\bf a}}_{u_1} & \bar{{\bf a}}_{u_1}& \bar{{\bf a}}_{u_3}& \bar{{\bf a}}_{u_3}& \bar{{\bf a}}_{u_2}& \bar{{\bf a}}_{u_2}
\end{pmatrix}
\]
respectively if $i\equiv k\mod 6$ for some $1 \leq k \leq 6$. Take
$$J_i=\{(x_1, x_2, x_3)\in A_i\mid (\exists \,t\in \{1, 2, 3\})~ x_t=0 \}.$$
It is easy to verify that $J_i$ is an ideal of $A_i$ and factor semiring $A_i/J_i$ is generated by $\{\bar{\alpha}_1, \ldots, \bar{\alpha}_6, \bar{\beta}_{11}, \bar{\beta}_{12}, \bar{\beta}_{13},\ldots , \bar{\beta}_{(i-1)1}, \bar{\beta}_{(i-1)2}, \bar{\beta}_{(i-1)3}\} \cup \{\bar{\beta}_{i1}, \bar{\beta}_{i2}, \bar{\beta}_{i3}\}$.
Also, it is easy to verify that $A_i/J_i$ is isomorphic to hypergraph semiring $S_{\overline{\scriptscriptstyle\mathbb{B}}_{i+1}}$, where hypergraph $\overline{\mathbb{B}}_{i+1}=(\overline{V}_{i+1}, \overline{E}_{i+1})$ with $\overline{V}_{i+1}=\overline{V}_{i}\cup \{\bar{\beta}_{i1}, \bar{\beta}_{i2}, \bar{\beta}_{i3}\}$
and
\begin{equation*}
\overline{E}_{i+1}= \begin{cases}
\overline{E}_{3}\cup \{\{\bar{\beta}_{22}, \bar{\beta}_{31}, \bar{\beta}_{32} \}, \{\bar{\beta}_{32}, \bar{\beta}_{33}, \bar{\alpha}_{3} \}\} & \text{if}~ i=3;\\
\overline{E}_{i}\cup \{\{ \bar{\beta}_{(i-1)2}, \bar{\beta}_{i1}, \bar{\beta}_{i2} \}, \{\bar{\beta}_{i2}, \bar{\beta}_{i3}, \bar{\beta}_{(i-2)2} \}\} & \text{if}~ i\geq 4 ~ \text{and} ~ i ~\text{is odd}~  ;\\
\overline{E}_{i}\cup \{\{\bar{\beta}_{(i-2)2}, \bar{\beta}_{i1},  \bar{\beta}_{i2} \}, \{\bar{\beta}_{i2}, \bar{\beta}_{i3}, \bar{\beta}_{(i-1)2} \}\}  & \text{if}~ i\geq 4 ~ \text{and} ~ i ~\text{is even},
\end{cases}
\end{equation*}
On the other hand, define the map $\varphi_{i+1}: \overline{V}_{i+1}\rightarrow {V}_{i+1}$ by
($\forall \,\, x\in \overline{V}_{i}\setminus \{\bar{\beta}_{11},\bar{\beta}_{12},\bar{\beta}_{13}\}$)\,\, $\varphi_{i+1}(x)= \varphi_{i}(x)$,
$\varphi_{i+1}(\bar{\beta}_{11})= u_{3i+4}, \,\, \varphi_{i+1}(\bar{\beta}_{12})= u_{3i+5}, \,\, \varphi_{i+1}(\bar{\beta}_{13})= u_{3i+6}$ and
$$\varphi_{i+1}(\bar{\beta}_{i1})= u_{3i+1}, \,\, \varphi_{i+1}(\bar{\beta}_{i2})= u_{3i+2}, \,\, \varphi_{i+1}(\bar{\beta}_{i3})= u_{3i+3}.$$
Then $\varphi_{i+1}$ is a hypergraph isomorphism  from $\overline{\mathbb{B}}_{i+1}$ to $\mathbb{B}_{i+1}$ and so $S_{\overline{\scriptscriptstyle\mathbb{B}}_{i+1}} \cong S_{{\scriptscriptstyle\mathbb{B}}_{i+1}}$.
This implies that   $A_{i}/J_{i}\cong S_{{\scriptscriptstyle\mathbb{B}}_{i+1}}$ and so $ S_{{\scriptscriptstyle\mathbb{B}}_{i+1}} \in  \mathbf{V}(S_{{\scriptscriptstyle \mathbb{B}}_{i}})$. This completes our proof.
\epf

Each 3-hypergraph $\mathbb{F}_i$ ($i \geq 1$) is given in following figure
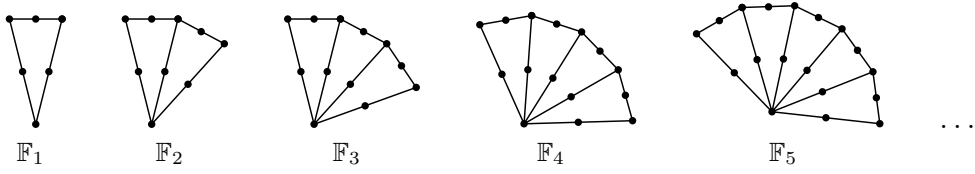
\begin{figure}[H]
	\centering
	\captionsetup[subfloat]{labelsep=none, format=plain, labelformat=empty}	
	\subfloat[$\mathbb{F}_1$]{
		\scalebox{0.7}{		
			\begin{tikzpicture}[ style=thick]
			\draw (0,0) -- (-0.5, 2) --(0.5,2) --cycle;
			\fill (0,0)  circle (2pt) ;
			\fill (-0.5, 2)  circle (2pt) ;
			\fill (0.5,2)  circle (2pt) ;
			\fill (0.5,2)   circle  (2pt) ;
			\fill (-0.25,1)  circle  (2pt) ;
			\fill (0.25,1)  circle  (2pt) ;
			\fill (0, 2)  circle  (2pt) ;
			\end{tikzpicture}
	}}
	\quad
	\subfloat[$\mathbb{F}_2$]{
		\scalebox{0.7}{	
			\begin{tikzpicture}[style=thick]
			\draw (0,0) -- (-0.5, 2) --(0.5,2) --cycle;
			\draw (0.5,2)--(1.38, 1.53);
			\draw (1.38, 1.53)--(0,0);
			\fill (0,0)  circle (2pt) ;
			\fill (-0.5, 2)  circle (2pt) ;
			\fill (0.5,2)  circle (2pt) ;
			\fill (0.5,2)   circle  (2pt) ;
			\fill (-0.25,1)  circle  (2pt) ;
			\fill (0.25,1)  circle  (2pt) ;
			\fill (0, 2)  circle  (2pt) ;
			\fill (1.38, 1.53)  circle  (2pt) ;
			\fill (0.94, 1.76)  circle  (2pt) ;
			\fill (0.69, 0.76)  circle  (2pt) ;
			\end{tikzpicture}
	}}	
	\quad
	\subfloat[$\mathbb{F}_3$]{
		\scalebox{0.7}{	
			\begin{tikzpicture}[style=thick]
			\draw (0,0) -- (-0.5, 2) --(0.5,2) --cycle;
			\draw (0.5,2)--(1.38, 1.53);
			\draw (1.38, 1.53)--(0,0);
			\draw (1.38, 1.53)--(1.94, 0.7);
			\draw (1.94, 0.7)--(0,0);
			\fill (0,0)  circle (2pt) ;
			\fill (-0.5, 2)  circle (2pt) ;
			\fill (0.5,2)  circle (2pt) ;
			\fill (0.5,2)   circle  (2pt) ;
			\fill (-0.25,1)  circle  (2pt) ;
			\fill (0.25,1)  circle  (2pt) ;
			\fill (0, 2)  circle  (2pt) ;
			\fill (1.38, 1.53)  circle  (2pt) ;
			\fill (0.94, 1.76)  circle  (2pt) ;
			\fill (0.69, 0.76)  circle  (2pt) ;
			\fill (1.94, 0.7)  circle  (2pt) ;
			\fill (1.66, 1.11)  circle  (2pt) ;
			\fill (0.97, 0.35)  circle  (2pt) ;
			\end{tikzpicture}
	}}	
	\quad
	\subfloat[$\mathbb{F}_4$]{
		\scalebox{0.7}{	
			\begin{tikzpicture}[rotate=10, style=thick]
			\draw (0,0) -- (-0.5, 2) --(0.5,2) --cycle;
			\draw (0.5,2)--(1.38, 1.53);
			\draw (1.38, 1.53)--(0,0);
			\draw (1.38, 1.53)--(1.94, 0.7);
			\draw (1.94, 0.7)--(0,0);
			\draw (1.94, 0.7)--(2.04,-0.3);
			\draw (2.04,-0.3)--(0,0);
			\fill (0,0)  circle (2pt) ;
			\fill (-0.5, 2)  circle (2pt) ;
			\fill (0.5,2)  circle (2pt) ;
			\fill (0.5,2)   circle  (2pt) ;
			\fill (-0.25,1)  circle  (2pt) ;
			\fill (0.25,1)  circle  (2pt) ;
			\fill (0, 2)  circle  (2pt) ;
			\fill (1.38, 1.53)  circle  (2pt) ;
			\fill (0.94, 1.76)  circle  (2pt) ;
			\fill (0.69, 0.76)  circle  (2pt) ;
			\fill (1.94, 0.7)  circle  (2pt) ;
			\fill (1.66, 1.11)  circle  (2pt) ;
			\fill (0.97, 0.35)  circle  (2pt) ;
			\fill (2.04,-0.3)  circle  (2pt) ;
			\fill (1.99, 0.2)  circle  (2pt) ;
			\fill (1.02, -0.15)  circle  (2pt) ;
			\end{tikzpicture}
	}}
	\quad
	\subfloat[$\mathbb{F}_5$]{
		\scalebox{0.7}{	
			\begin{tikzpicture}[rotate=30, style=thick]
			\draw (0,0) -- (-0.5, 2) --(0.5,2) --cycle;
			\draw (0.5,2)--(1.38, 1.53);
			\draw (1.38, 1.53)--(0,0);
			\draw (1.38, 1.53)--(1.94, 0.7);
			\draw (1.94, 0.7)--(0,0);
			\draw (1.94, 0.7)--(2.04,-0.3);
			\draw (2.04,-0.3)--(0,0);
			\draw (2.04,-0.3)--(1.66,-1.22);
			\draw (1.66,-1.22)--(0,0);
			\fill (0,0)  circle (2pt) ;
			\fill (-0.5, 2)  circle (2pt) ;
			\fill (0.5,2)  circle (2pt) ;
			\fill (0.5,2)   circle  (2pt) ;
			\fill (-0.25,1)  circle  (2pt) ;
			\fill (0.25,1)  circle  (2pt) ;
			\fill (0, 2)  circle  (2pt) ;
			\fill (1.38, 1.53)  circle  (2pt) ;
			\fill (0.94, 1.76)  circle  (2pt) ;
			\fill (0.69, 0.76)  circle  (2pt) ;
			\fill (1.94, 0.7)  circle  (2pt) ;
			\fill (1.66, 1.11)  circle  (2pt) ;
			\fill (0.97, 0.35)  circle  (2pt) ;
			\fill (2.04,-0.3)  circle  (2pt) ;
			\fill (1.99, 0.2)  circle  (2pt) ;
			\fill (1.02, -0.15)  circle  (2pt) ;
			\fill (1.66,-1.22)  circle  (2pt) ;
			\fill (0.83, -0.61)  circle  (2pt) ;
			\fill (1.85, -0.76)  circle  (2pt) ;
			\end{tikzpicture}
	}}
	\hspace{5mm} \ldots \hspace{5mm}
	\caption{Fan-type 3-hypergraph}\label{tu2}
\end{figure}
\vspace{-3mm}
\noindent where all hyperedges in each $\mathbb{F}_i$ are only denoted by straight lines and all vertices are denoted by black dots. We refer to $\mathbb{F}_i$ as the $i$-th fan-type 3-hypergraph.

Similarly to Proposition \ref{prop3.4}, we immediately have

\begin{prop}
	For all $i\geq 1$,
	$\mathbf{V}(S_{{\scriptscriptstyle \mathbb{F}}_{i}})= \mathbf{V}(S_{{\scriptscriptstyle \mathbb{B}}_{1}})$.
\end{prop}

Each 3-hypergraph $\mathbb{N}_i$ ($i \geq 1$) is given in following figure
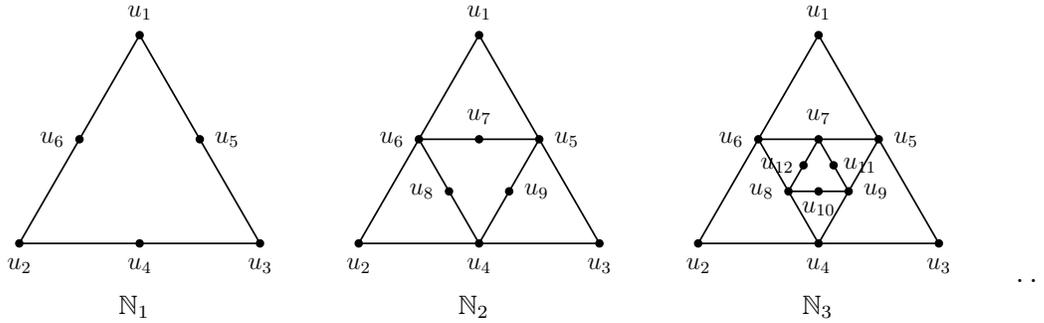
\begin{figure}[H]
	\centering
	\captionsetup[subfloat]{labelsep=none, format=plain, labelformat=empty}	
	\subfloat[$\mathbb{N}_1$]{  	
		\scalebox{0.8}{	
			\begin{tikzpicture}[style=thick]
			\draw (0,0) -- (1,0) -- (2,0) -- (3,0) -- (4,0) -- (3, 1.732050807568877) -- (2, 3.464101615137754) -- (1, 1.732050807568877) --cycle;	
			\fill (2, 3.464101615137754) node[above=1mm]{$u_1$} circle (2pt) ;
			\fill (1,1.732050807568877) node[left=1mm]{$u_6$} circle  (2pt) ;
			\fill (0,0) node[below=1mm]{$u_2$} circle (2pt) ;
			\fill (2,0) node[below=1mm]{$u_4$} circle (2pt) ;
			\fill (4,0) node[below=1mm]{$u_3$} circle (2pt) ;
			\fill (3,1.732050807568877) node[right=1mm]{$u_5$} circle  (2pt) ;
			\end{tikzpicture}
	}}
	\quad
	\subfloat[$\mathbb{N}_2$]{
		\scalebox{0.8}{		
			\begin{tikzpicture}[style=thick]
			\draw (0,0) -- (1,0) -- (2,0) -- (3,0) -- (4,0) -- (3, 1.732050807568877) -- (2, 3.464101615137754) -- (1, 1.732050807568877) --cycle;
			\draw (2,0) -- (1.5,0.8660254037844386) -- (1,1.732050807568877) -- (2, 1.732050807568877)-- (3, 1.732050807568877) --  (2.5,0.8660254037844386) -- cycle;		
			\fill (2, 3.464101615137754) node[above=1mm]{$u_1$} circle (2pt) ;
			\fill (1,1.732050807568877) node[left=1mm]{$u_6$} circle  (2pt) ;
			\fill (0,0) node[below=1mm]{$u_2$} circle (2pt) ;
			\fill (2,0) node[below=1mm]{$u_4$} circle (2pt) ;
			\fill (4,0) node[below=1mm]{$u_3$} circle (2pt) ;
			\fill (3,1.732050807568877) node[right=1mm]{$u_5$} circle  (2pt) ;
			\fill (2, 1.732050807568877) node[above=1mm]{$u_7$} circle (2pt) ;
			\fill (1.5,0.8660254037844386) node[left=1mm]{$u_8$} circle  (2pt) ;
			\fill (2.5,0.8660254037844386) node[right=1mm]{$u_9$} circle  (2pt) ;
			\end{tikzpicture}	}}	
	\quad
	\subfloat[ $\mathbb{N}_3$]{
		\scalebox{0.8}{		
			\begin{tikzpicture}[style=thick]
			\draw (0,0) -- (1,0) -- (2,0) -- (3,0) -- (4,0) -- (3, 1.732050807568877) -- (2, 3.464101615137754) -- (1, 1.732050807568877) --cycle;
			\draw (2,0) -- (1.5,0.8660254037844386) -- (1,1.732050807568877) -- (2, 1.732050807568877)-- (3, 1.732050807568877) --  (2.5,0.8660254037844386) -- cycle;		
			\draw (1.5,0.8660254037844386) -- (1.75,1.299038105676658) -- (2, 1.732050807568877)-- (2.25, 1.299038105676658) --  (2.5,0.8660254037844386) --  (2,0.8660254037844386)-- cycle;		
			\fill (2, 3.464101615137754) node[above=1mm]{$u_1$} circle (2pt) ;
			\fill (1,1.732050807568877) node[left=1mm]{$u_6$} circle  (2pt) ;
			\fill (0,0) node[below=1mm]{$u_2$} circle (2pt) ;
			\fill (2,0) node[below=1mm]{$u_4$} circle (2pt) ;
			\fill (4,0) node[below=1mm]{$u_3$} circle (2pt) ;
			\fill (3,1.732050807568877) node[right=1mm]{$u_5$} circle  (2pt) ;
			\fill (2, 1.732050807568877) node[above=1mm]{$u_7$} circle (2pt) ;
			\fill (1.5,0.8660254037844386) node[left=1mm]{$u_8$} circle  (2pt) ;
			\fill (2.5,0.8660254037844386) node[right=1mm]{$u_9$} circle  (2pt) ;
			\fill (1.75,1.299038105676658) node[left=0.1mm]{$u_{12}$} circle  (2pt) ;
			\fill (2,0.8660254037844386) node[below=0.1mm]{$u_{10}$} circle  (2pt) ;
			\fill (2.25,1.299038105676658) node[right=0.1mm]{$u_{11}$} circle  (2pt) ; 	
			\end{tikzpicture}
	}}
	\hspace{5mm} \ldots \hspace{5mm}
	\caption{Nested-type 3-hypergraph}\label{tu3}
\end{figure}
\vspace{-3mm}
\noindent where all hyperedges in each $\mathbb{N}_i$ are only denoted by straight lines and all vertices are denoted by black dots. We refer to $\mathbb{N}_i$ as the $i$-th nested-type 3-hypergraph, and refer to the subhypergraph of $\mathbb{N}_{i}$ induced by vertex set $\{u_{3i+3},  u_{3i+2}, u_{3i+1}, u_{3i}, u_{3i-1}, u_{3i-2}\}$ as the $i$-th triangle. Notice that the hyperedge set of the $i$-th triangle is
$$\{\{u_{3i-1}, u_{3i+3}, u_{3i-2}\}, \{u_{3i-2}, u_{3i+2}, u_{3i}\}, \{u_{3i}, u_{3i+1}, u_{3i-1}\} \} .$$

In the remainder of this section we shall study the varieties $\mathbf{V} (S_{\scriptscriptstyle \mathbb{N}_i})$. The following Proposition \ref{xinprop4.5} will tell us that
$$\mathbf{V} (S_{\scriptscriptstyle \mathbb{B}_1}) = \mathbf{V} (S_{\scriptscriptstyle \mathbb{N}_1}) \subset \mathbf{V} (S_{\scriptscriptstyle \mathbb{N}_2}) \subset \cdots \subset \mathbf{V} (S_{\scriptscriptstyle \mathbb{N}_i})\subset \cdots.$$
That is to say, it is
an infinite ascending chain in the lattice of subvarieties of the variety generated by all 3-hypergraph semirings.
To prove that $\mathbf{V}(S_{{\scriptscriptstyle \mathbb{N}}_{i}})\subset \mathbf{V}(S_{{\scriptscriptstyle \mathbb{N}}_{i+1}})$ for all $i\geq 1$, the following identities are needed for us
\begin{align}
p_1 &\approx q_1  \label{gongshi4.1} \\
p_1+p_2 &\approx q_1+q_2 \label{gongshi4.2} \\
p_1+p_2+p_3 &\approx q_1+q_2+q_3 \label{gongshi4.3}\\
&~~  \vdots \notag
\end{align}
where
\begin{equation*}
\begin{split}
p_1 &= x_1x_6x_2+x_2x_4x_3+x_3x_5x_1 \\
p_2 &= x_4x_9x_5+x_5x_7x_6+x_6x_8x_4 \\
p_3 &= x_7x_{12}x_8+x_8x_{10}x_9+x_9x_{11}x_7\\
&~~~~~~~~~~~~~~~~~ \vdots \\
q_1 &= (x_1+x_4)(x_2+x_5)(x_3+x_6)\\
q_2 &= (x_4+x_7)(x_5+x_8)(x_6+x_9)\\
q_3 &= (x_7+x_{10})(x_8+x_{11})(x_9+x_{12})\\
&~~~~~~~~~~~~~~~~~  \vdots
\end{split}
\end{equation*}

Note that in the above the $i$-th identity is $\sum_{j=1}^i p_j \approx \sum_{j=1}^i q_j $, where $p_i=x_{3i-1} x_{3i+3} x_{3i-2}+ x_{3i-2} x_{3i+2} x_{3i} +x_{3i} x_{3i+1} x_{3i-1}$ and
$q_i=(x_{3i-1}+x_{3i+2})(x_{3i-2}+x_{3i+1})(x_{3i}+x_{3i+3})$.

\begin{prop} \label{xinprop4.5}
	For all $i\geq 1$,	
	$\mathbf{V}(S_{{\scriptscriptstyle \mathbb{N}}_{i}})\subset \mathbf{V}(S_{{\scriptscriptstyle \mathbb{N}}_{i+1}})$.
\end{prop}
\pf
It is easy to see that $\mathbf{V}(S_{{\scriptscriptstyle \mathbb{N}}_{i}})\subseteq \mathbf{V}(S_{{\scriptscriptstyle \mathbb{N}}_{i+1}})$ for all $i\geq 1$.
We need only to prove that $\mathbf{V}(S_{{\scriptscriptstyle \mathbb{N}}_{i}})\neq \mathbf{V}(S_{{\scriptscriptstyle \mathbb{N}}_{i+1}})$.

(1) If $i=1$, then $S_{{\scriptscriptstyle \mathbb{N}}_1}$ satisfies the identity $(\ref{gongshi4.2})$
but $S_{{\scriptscriptstyle \mathbb{N}}_2}$ does not.
This  means that $\mathbf{V}(S_{{\scriptscriptstyle \mathbb{N}}_{1}})\neq \mathbf{V}(S_{{\scriptscriptstyle \mathbb{N}}_{2}})$. In fact, it is easy to verify that $S_{{\scriptscriptstyle \mathbb{N}}_2}$ does not satisfy the identity $(\ref{gongshi4.2})$.
To show that  $S_{{\scriptscriptstyle \mathbb{N}}_1}$ satisfies the identity $(\ref{gongshi4.2})$,
we need only to prove that for any given substitution $\varphi$  from $\{x_1, x_2, \ldots,x_9 \}$ to $S_{{\scriptscriptstyle \mathbb{N}}_1}$,
$$\varphi(p_1+p_2)=\bar{{\bf a}}\iff \varphi(q_1+q_2)=\bar{{\bf a}},$$
since
$\varphi(p_1+p_2), \varphi(q_1+q_2)\in \{0, \bar{{\bf a}} \}.$

(i) If $\varphi(p_1+p_2)=\bar{{\bf a}}$, then $\varphi(x_rx_sx_t)=\bar{{\bf a}}$ for each word $x_r x_s x_t$ occurring in $p_1$ and $p_2$.
In particular, $\varphi(x_1x_6x_2)=\bar{{\bf a}}$.
Without loss of generality, assume that $\{\varphi(x_1), \varphi(x_6), \varphi(x_2)\} = \{\bar{{\bf a}}_{u_1},  \bar{{\bf a}}_{u_6},  \bar{{\bf a}}_{u_2}\}$, since $\bar{{\bf a}} = \bar{{\bf a}}_{u_1} \bar{{\bf a}}_{u_6} \bar{{\bf a}}_{u_2}$.

It is enough to consider only the following two cases.

{\bfseries Case 1.} If $\varphi(x_1)=\bar{{\bf a}}_{u_1}$, $\varphi(x_6)=\bar{{\bf a}}_{u_6}$, $\varphi(x_2)=\bar{{\bf a}}_{u_2}$,
then $\varphi(x_5)\varphi(x_3)=\bar{{\bf a}}_{u_2} \bar{{\bf a}}_{u_6}$, $\varphi(x_4)\varphi(x_3)=\bar{{\bf a}}_{u_1} \bar{{\bf a}}_{u_6}$, $\varphi(x_8)\varphi(x_4)=\varphi(x_7)\varphi(x_5)=\bar{{\bf a}}_{u_1} \bar{{\bf a}}_{u_2}$.
This implies that  $$\varphi(x_1)=\varphi(x_4)=\varphi(x_7)=\bar{{\bf a}}_{u_1},\,\, \varphi(x_6)=\varphi(x_3)=\varphi(x_9)=\bar{{\bf a}}_{u_6} ,$$
and $\varphi(x_2)=\varphi(x_5)=\varphi(x_8)=\bar{{\bf a}}_{u_2}.$
Thus $\varphi(q_1+q_2)=\bar{{\bf a}}$, as required.

{\bfseries Case 2.} If $\varphi(x_1)=\bar{{\bf a}}_{u_1}$, $\varphi(x_6)=\bar{{\bf a}}_{u_2}$, $\varphi(x_2)=\bar{{\bf a}}_{u_6}$,
then $\varphi(x_5)\varphi(x_3)=\bar{{\bf a}}_{u_2} \bar{{\bf a}}_{u_6}$, $\varphi(x_4)\varphi(x_3)=\bar{{\bf a}}_{u_1} \bar{{\bf a}}_{u_2}$, $\varphi(x_8)\varphi(x_4)=\varphi(x_7)\varphi(x_5)=\bar{{\bf a}}_{u_1} \bar{{\bf a}}_{u_6}$.
This implies that  $$\varphi(x_1)=\varphi(x_4)=\varphi(x_7)=\bar{{\bf a}}_{u_1}, \,\,\varphi(x_6)=\varphi(x_3)=\varphi(x_9)=\bar{{\bf a}}_{u_2},$$
and $\varphi(x_2)=\varphi(x_5)=\varphi(x_8)=\bar{{\bf a}}_{u_6}$.
Thus  $\varphi(q_1+q_2)=\bar{{\bf a}}$, as required.

(ii) If $\varphi(q_1+q_2)=\bar{{\bf a}}$, then  $$\varphi(x_1)=\varphi(x_4)=\varphi(x_7),\,\, \varphi(x_6)=\varphi(x_3)=\varphi(x_9),\,\, \varphi(x_2)=\varphi(x_5)=\varphi(x_8), $$
and $\varphi(x_1x_6x_2)=\bar{{\bf a}}$. This implies that  $\varphi(p_1 + p_2)=\bar{{\bf a}}$, as required.

Summarizing (i) and (ii) we have that $\varphi(p_1+p_2)=\bar{{\bf a}}$ if and only if $\varphi(q_1+q_2)=\bar{{\bf a}}$.
This shows that  $\mathbf{V}(S_{{\scriptscriptstyle \mathbb{N}}_{1}})\subset \mathbf{V}(S_{{\scriptscriptstyle \mathbb{N}}_{2}})$.

(2) If $i=2$, then $S_{{\scriptscriptstyle \mathbb{N}}_2}$ satisfies the identity  $(\ref{gongshi4.3})$ but $S_{{\scriptscriptstyle \mathbb{N}}_3}$ does not.
This means that $\mathbf{V}(S_{{\scriptscriptstyle \mathbb{N}}_{2}})\neq \mathbf{V}(S_{{\scriptscriptstyle \mathbb{N}}_{3}})$.
In fact, it is easy to verify that $S_{{\scriptscriptstyle \mathbb{N}}_3}$ does not satisfy the identity $(\ref{gongshi4.3})$.
To show that $S_{{\scriptscriptstyle \mathbb{N}}_2}$ satisfies the identity  $(\ref{gongshi4.3})$, we need only to prove that  for any given substitution $\varphi$ from $\{x_1, x_2, \ldots,x_{12} \}$ to $S_{{\scriptscriptstyle \mathbb{N}}_2}$,
$$\varphi(p_1+p_2+p_3)=\bar{{\bf a}} \iff \varphi(q_1+q_2+q_3)=\bar{{\bf a}},$$
since
$\varphi(p_1+p_2+p_3), \varphi(q_1+q_2+q_3)\in \{0, \bar{{\bf a}} \}$.

(i) If $\varphi(p_1+p_2+p_3)=\bar{{\bf a}}$, then $\varphi(x_rx_sx_t)=\bar{{\bf a}}$ for each word $x_r x_s x_t$ occurring in $p_1$, $p_2$ and $p_3$.
In particular, $\varphi(x_1x_6x_2)=\bar{{\bf a}}$.
Without loss of generality, assume that $\{\varphi(x_1), \varphi(x_6), \varphi(x_2)\} = \{\bar{{\bf a}}_{u_1},  \bar{{\bf a}}_{u_6},  \bar{{\bf a}}_{u_2}\}$, since $\bar{{\bf a}} = \bar{{\bf a}}_{u_1} \bar{{\bf a}}_{u_6} \bar{{\bf a}}_{u_2}$.

It is enough to consider only the following two cases.

{\bfseries Case 1.} If $\varphi(x_1)=\bar{{\bf a}}_{u_1}$, $\varphi(x_6)=\bar{{\bf a}}_{u_6}$, $\varphi(x_2)=\bar{{\bf a}}_{u_2}$,
then $\varphi(x_5)\varphi(x_3)=\bar{{\bf a}}_{u_2} \bar{{\bf a}}_{u_6}$, $\varphi(x_4)\varphi(x_3)=\bar{{\bf a}}_{u_1} \bar{{\bf a}}_{u_6}$, $\varphi(x_8)\varphi(x_4)=\varphi(x_7)\varphi(x_5)=\bar{{\bf a}}_{u_1} \bar{{\bf a}}_{u_2}$.
This implies that  $$\varphi(x_1)=\varphi(x_4)=\varphi(x_7)=\varphi(x_{10})=\bar{{\bf a}}_{u_1} ,\,\, \varphi(x_6)=\varphi(x_3)=\varphi(x_9)=\varphi(x_{12})=\bar{{\bf a}}_{u_6}, $$
and
$\varphi(x_2)=\varphi(x_5)=\varphi(x_8)=\varphi(x_{11})=\bar{{\bf a}}_{u_2}$.
Thus $\varphi(q_1+q_2+q_3)=\bar{{\bf a}}$, as required.

{\bfseries Case 2.} If $\varphi(x_1)=\bar{{\bf a}}_{u_1}$, $\varphi(x_6)=\bar{{\bf a}}_{u_2}$, $\varphi(x_2)=\bar{{\bf a}}_{u_6}$,
then $\varphi(x_5)\varphi(x_3)=\bar{{\bf a}}_{u_2} \bar{{\bf a}}_{u_6}$, $\varphi(x_4)\varphi(x_3)=\bar{{\bf a}}_{u_1} \bar{{\bf a}}_{u_2}$, $\varphi(x_8)\varphi(x_4)=\varphi(x_7)\varphi(x_5)=\bar{{\bf a}}_{u_1} \bar{{\bf a}}_{u_6}$.
This implies that  $$\varphi(x_1)=\varphi(x_4)=\varphi(x_7)=\varphi(x_{10})=\bar{{\bf a}}_{u_1} , \,\, \varphi(x_6)=\varphi(x_3)=\varphi(x_9)=\varphi(x_{12})=\bar{{\bf a}}_{u_2} ,$$
and $\varphi(x_2)=\varphi(x_5)=\varphi(x_8)=\varphi(x_{11})=\bar{{\bf a}}_{u_6}$.
Thus $\varphi(q_1+q_2+q_3)=\bar{{\bf a}}$, as required.

(ii) If $\varphi(q_1+q_2+q_3)=\bar{{\bf a}}$, then $\varphi(x_1)=\varphi(x_4)=\varphi(x_7)=\varphi(x_{10})$,
$$\varphi(x_6)=\varphi(x_3)=\varphi(x_9)=\varphi(x_{12}), \,\, \varphi(x_2)=\varphi(x_5)=\varphi(x_8)=\varphi(x_{11}) , $$
and $\varphi(x_1x_6x_2)=\bar{{\bf a}}$. This implies that  $\varphi(p_1+p_2+p_3)=\bar{{\bf a}}$, as required.

Summarizing (i) and (ii) we have that $\varphi(p_1+p_2+p_3)=\bar{{\bf a}}\,$ if and only if $\,\varphi(q_1+q_2+q_3)=\bar{{\bf a}}$. This shows that  $\mathbf{V}(S_{{\scriptscriptstyle \mathbb{N}}_{2}})\subset \mathbf{V}(S_{{\scriptscriptstyle \mathbb{N}}_{3}})$.

Similarly, we can prove that
$S_{{\scriptscriptstyle \mathbb{N}}_i}$ satisfies the identity
$\sum_{j=1}^{i+1} p_j \approx \sum_{j=1}^{i+1} q_j $ but  $S_{{\scriptscriptstyle \mathbb{N}}_{i+1}}$ does not.
This shows that
$\mathbf{V}(S_{{\scriptscriptstyle \mathbb{N}}_{i}})\subset \mathbf{V}(S_{{\scriptscriptstyle \mathbb{N}}_{i+1}})$ for all $i \geq 1$. This completes our proof.
\epf
\noindent {\bf Remark} The above proposition tells us that the variety generated by all 3-hypergraph semirings with girth being equal to 3 has infinitely many subvarieties.

From Proposition \ref{lem4}, we have that $S_{{\scriptscriptstyle \mathbb{N}}_{1}}\in \mathbf{V}(S_c(abcd))$. However, the following consequences show that $S_{{\scriptscriptstyle \mathbb{N}}_{i}}\notin \mathbf{V}(S_c(abcd))$ for all $i\geq 2$.
Recall from \cite{Jackson} that $S_7$ is a 3-element ai-semiring $\{1, a, 0\}$ with the Cayley tables
\begin{center}
\begin{tabular}{c|ccc}
			$+$&1&$a$&0\\
			\hline
			$1$&1&$0$&0\\
			$a$&0&$a$&0\\
			$0$&0&$0$&0\\
\end{tabular}\qquad
\begin{tabular}{c|ccc}
			$\cdot$&1&$a$&0\\
			\hline
			$1$&1&$a$&0\\
			$a$&$a$&$0$&0\\
			$0$&0&$0$&0\\
\end{tabular}~.
\end{center}

\begin{prop}\label{prop4.6}
	For all $i\geq 2$, neither $S_{{\scriptscriptstyle \mathbb{N}}_{i}}$ nor  $S_7$ lies in the variety generated by the other.
\end{prop}
\pf
It is easy to verify that  $S_7$ satisfies the identity $(\ref{gongshi4.2})$, while $S_{{\scriptscriptstyle \mathbb{N}}_{i}}$ does not. Thus $S_{{\scriptscriptstyle \mathbb{N}}_{i}} \notin \mathbf{V}(S_7)$.
Conversely, $S_{{\scriptscriptstyle \mathbb{N}}_{i}}$ satisfies the following identity  $(\ref{gongshi4.4})$ but $S_7$ does not,
\begin{equation}\label{gongshi4.4}
x_1x_2x_3x_4 \approx y_1y_2y_3y_4.
\end{equation} 
This implies that $S_7 \notin \mathbf{V}(S_{{\scriptscriptstyle \mathbb{N}}_{i}})$.
\epf

From \cite[Lemma 2.10]{Jackson}, we know that $S_c(abcd)\in  \mathbf{V}(S_7)$. Furthermore, by above proposition, we have that $S_{{\scriptscriptstyle \mathbb{N}}_{i}} \notin \mathbf{V}(S_c(abcd))$ for all $i\geq 2$.  
It is also easy to see that  $S_{{\scriptscriptstyle \mathbb{N}}_i}$ satisfies the following identity  $(\ref{gongshi4.4}) $ but $S_c(abcd)$ does not. This implies that $S_c(abcd) \notin \mathbf{V}(S_{{\scriptscriptstyle \mathbb{N}}_{i}})$. 
Thus we have the following corollary to Proposition \ref{prop4.6}.

\begin{cor}
	For all $i\geq 2$, neither $S_{{\scriptscriptstyle \mathbb{N}}_{i}}$ nor  $S_c(abcd)$ lies in the variety generated by the other.
\end{cor}

\noindent {\bf Remark} We feel that it is interesting to study the lattice of subvarieties of variety generated by all 3-hypergraph semirings and hard to give a characterization of it. This paper is only first step in the topic. In fact, in the next manuscript (as a sisters of this paper), we will give an infinite descending chain of varieties generated by 3-hypergraph semirings in the lattice, which also is in interval $[\mathbf{V}(S_c(abc)),\,\mathbf{V}(S_c(abcd))]$. This shows that
different from $\mathbf{V}(S_c(abc))$ (see, \cite{rz}), the lattice of subvarieties of $\mathbf{V}(S_c(abcd))$ is infinite.

\end{document}